\documentclass[11pt]{amsart}
\usepackage[T1]{fontenc}
\usepackage{lmodern}
\usepackage{amsmath,amssymb,amsthm,mathtools}
\usepackage[margin=1.08in]{geometry}
\usepackage{microtype,enumitem,booktabs,needspace}
\usepackage[hidelinks]{hyperref}
\hypersetup{pdftitle={Rigidity of free boundary minimal Y-surfaces in the unit ball},
 pdfauthor={Elham Matinpour},pdfsubject={Conformal, geometric, and Morse index rigidity},
 pdfkeywords={free boundary minimal surfaces, Y-singularities, rigidity, Morse index, nullity, Steklov spectrum}}
\newtheorem{theorem}{Theorem}[section]
\newtheorem{proposition}[theorem]{Proposition}
\newtheorem{lemma}[theorem]{Lemma}
\newtheorem{corollary}[theorem]{Corollary}
\theoremstyle{definition}
\newtheorem{definition}[theorem]{Definition}
\newtheorem{question}[theorem]{Question}
\theoremstyle{remark}
\newtheorem{remark}[theorem]{Remark}
\numberwithin{equation}{section}
\newcommand{\R}{\mathbb{R}}\newcommand{\C}{\mathbb{C}}
\newcommand{\B}{\mathbb{B}}\newcommand{\D}{\mathbb{D}}
\newcommand{\Sph}{\mathbb{S}}\newcommand{\II}{\mathrm{II}}
\newcommand{\ip}[2]{\langle #1,#2\rangle}
\newcommand{\Hopf}{\mathcal H}\newcommand{\sigmasp}{\widehat\sigma}
\newcommand{\nsc}{N_{\mathrm{sc}}}
\newcommand{\V}{\mathcal V}\newcommand{\Aop}{\mathcal A}
\DeclareMathOperator{\Area}{Area}\DeclareMathOperator{\Length}{Length}
\DeclareMathOperator{\Rea}{Re}\DeclareMathOperator{\Ima}{Im}
\DeclareMathOperator{\spanop}{span}\DeclareMathOperator{\Ind}{ind}
\DeclareMathOperator{\nul}{nul}\DeclareMathOperator{\diver}{div}

\newcommand{\proofstep}[2]{\par\smallskip\noindent Step~#1. \emph{#2}\ }
\title[Rigidity of free boundary minimal $Y$-surfaces]
{Rigidity of free boundary minimal $Y$-surfaces in the unit ball}
\author{Elham Matinpour}
\address{Department of Mathematics, Johns Hopkins University, 3400 N. Charles Street, Baltimore, MD 21218}
\email{ematinp1@jhu.edu}
\date{}
\subjclass[2020]{Primary 53A10; Secondary 49Q05, 49R05}
\keywords{Free boundary minimal surfaces, $Y$-singularities, rigidity, Morse index, nullity, Steklov spectrum}
\begin{document}
\begin{abstract}
We study rigidity for free boundary minimal $Y$-surfaces in the unit
ball $\B^n$, $n\geq3$, whose three disk faces meet at $120^\circ$
along an arc. The face immersions are $C^{2,\alpha}$ and remain of
rank two at the two junction endpoints. For a conformal minimal
immersion of the standard model of three half-disks, matching
parametrizations on the common diameter force the image to be the
flat $Y$-cone. We then give geometric criteria that do not require
this matching. If $c_j=\ip{\II_j(T,T)}{\II_j(T,\tau_j)}$ along the
junction, rigidity holds when one $c_j$ changes sign at most twice,
or when all three do so at most three times. In $\B^3$, a single sign
change of the common geodesic torsion is still enough. For the area
problem in which the junction may move and the outer boundary remains
on the sphere, every such surface has Morse index at least $n-1$, and
equality characterizes the flat $Y$-cone. Its normal nullity is
$3(n-2)$, and the null fields are precisely the normal parts of
ambient rotations. The arguments combine holomorphic differentials
and boundary moment identities with a coupled Steklov calculation and
the second variation of area.
\end{abstract}
\maketitle

\section{Introduction}
A free boundary minimal surface in the unit ball $\B^n\subset\R^n$
has boundary on $\Sph^{n-1}$ and meets the sphere orthogonally.
Nitsche~\cite{nitsche1985stationary} proved that such a minimal disk in
$\B^3$ is an equatorial disk. Souam~\cite{souam1997stability} obtained
related results in three-dimensional space forms. Fraser and
Schoen~\cite{fraser2015uniqueness} proved disk rigidity in balls in
space forms of arbitrary dimension.

We study three minimal disk faces meeting along a common arc $\Gamma$,
with outward unit conormals satisfying
\begin{equation}\label{eq:intro-balance}
 \tau_1+\tau_2+\tau_3=0.
\end{equation}
For unit conormals this is the $120^\circ$ condition for equal tensions.
The other edge of each face lies on the unit sphere and satisfies the
free boundary condition. We assume that the immersions, boundary arcs,
and junction identifications are $C^{2,\alpha}$, and that the
immersions remain of rank two at the vertices in ordinary extendible
corner coordinates. Immersed images are allowed, and area is counted
on the abstract faces. Definition~\ref{def:freeboundary} gives the
precise class, which we call regular disk-type free boundary minimal
$Y$-surfaces. Its regularity is part of the hypotheses rather than a
consequence of weak stationarity.

The standard compact $Y$-cone, denoted by $YC$, consists of three
half-unit disks joined along a diameter at $120^\circ$. Its linear
span is a copy of $\R^3$ in $\R^n$. We give three different sufficient
conditions for a surface in the stated class to equal this model.

Related rigidity results concern other boundary and completeness
conditions. Bernstein and Maggi~\cite{bernstein2021plateau} study
minimal surfaces with Plateau-like singularities, including symmetric
boundary configurations and complete surfaces with two regular ends.
Wang~\cite{wang2025bernstein} proves flatness for complete orientable
stable minimal Plateau surfaces with quadratic area growth in $\R^3$.
Here the surfaces are compact, and their outer boundary is constrained
to a fixed sphere.

\subsection{Conformal rigidity}
The first theorem concerns a conformal immersion of the standard
abstract model. Its three face maps must agree at corresponding points
of the identified diameter.

\begin{theorem}\label{thm:main1}
Let $n\geq3$ and $0<\alpha<1$. Suppose
$u=(u_1,u_2,u_3):YC\to\overline{\B^n}$ has face restrictions that are
$C^{2,\alpha}$ on the closed half-disks, conformal and minimal in their
interiors, and of rank two everywhere, including at the corners.
Assume that the maps agree pointwise on the identified diameter,
satisfy \eqref{eq:intro-balance}, and have spherical free boundary on
their outer semicircles. Then $u(YC)$ is the standard flat $Y$-cone,
up to an orthogonal transformation of $\R^n$.
\end{theorem}

Conformality on the individual faces is not enough for this theorem:
separate conformal parametrizations may induce different parameters
on the same junction. We require parametrizations that agree there.
Mese and Yamada~\cite{mese2006parameterized,mese2011local} study
conformal structures and local parametrizations for singular surfaces,
but their local result~\cite[Theorem~6]{mese2011local} does not
provide the global matching used here.
Wang~\cite{wang2021uniformization} similarly distinguishes
uniformization with fixed junction identifications from weak
uniformization, where those identifications may change.

In separate work~\cite{matinpour2026plateau}, the author studies the
tetrahedral $T$-cone and uses the conformal $Y$-model theorem.
The present paper retains that theorem and develops additional
geometric and variational criteria for the three-face junction.

\subsection{Geometric junction criteria}
Let $T$ be a consistently oriented unit tangent to $\Gamma$, and let
$\II_j$ be the normal-valued second fundamental form of face $j$.
Set
\begin{equation}\label{eq:intro-c}
 c_j=\ip{\II_j(T,T)}{\II_j(T,\tau_j)}.
\end{equation}
We use $\nsc$ for the number of strict alternations of nonzero signs
along the junction, as defined in Definition~\ref{def:signchanges}.
A zero without a change of sign does not count. Reversing $T$ reverses
the junction order and changes each $c_j$ to $-c_j$, leaving its count
unchanged.

\begin{theorem}\label{thm:oscillation}
Let $n\geq3$, and let $\Sigma\subset\overline{\B^n}$ satisfy
Definition~\ref{def:freeboundary}. With $c_j$ as in
\eqref{eq:intro-c}, either of the following conditions implies that
$\Sigma$ is the standard flat $Y$-cone, up to an orthogonal transformation.
\begin{enumerate}[label=\textup{(\roman*)},leftmargin=2.4em]
 \item $\nsc(c_j)\leq2$ for at least one face.
 \item $\nsc(c_j)\leq3$ for all three faces.
\end{enumerate}
\end{theorem}
The conformal face parametrizations in this theorem are chosen
separately. They are not required to agree along the junction.

\begin{corollary}\label{cor:geometric}
In the setting of Theorem~\ref{thm:oscillation}, $\Sigma$ is the
standard flat $Y$-cone if one face is planar, or if, on one face,
\[
 \II_j(T,\tau_j)\equiv0
 \quad\text{or}\quad
 \ip{\II_j(T,T)}{\II_j(T,\tau_j)}\equiv0
 \qquad\text{on }\Gamma.
\]
Consequently, on a nonflat surface satisfying
Definition~\ref{def:freeboundary}, every $c_j$ has at least three sign changes and at least
one has four or more.
\end{corollary}

In $\R^3$, orient the face normals along the junction by
$N_j=T\times\tau_j$. The quantities $\ip{\II_j(T,\tau_j)}{N_j}$
agree for the three faces, and we denote their common value by $b$.

\begin{theorem}\label{thm:torsion}
Let $\Sigma\subset\overline{\B^3}$ satisfy
Definition~\ref{def:freeboundary}. If $\nsc(b)\leq1$, then $\Sigma$
is the standard flat $Y$-cone.
\end{theorem}
This theorem also allows different conformal parameters on the three
junction edges. The required separate conformal charts follow from
the geometric regularity and the right angles at the vertices.
Appendix~\ref{app:corners} proves their boundary regularity.

Several ingredients in these boundary arguments have classical
precursors. Brander and L\'opez~\cite{brander2017remarks} derive
weighted identities for geodesic torsion and normal curvature on
constant mean curvature disks with analytic boundary. Choe's
theorem~\cite{choe2002round} gives umbilicity for an immersed constant
mean curvature disk in $\R^3$ when the regular boundary arcs are
curvature lines and there are at most three vertices of angle less
than $\pi$. In our three-dimensional setting that theorem applies
facewise when the junction is itself a curvature line, because the
spherical free boundary is then another curvature line. Here the
junction torsion may be nonzero, so we instead compare endpoint
identities from the three separate half-disks through the geometric
balance conditions. The moment argument used to control sign changes
is classical in spirit. For related background, see the
Sturm--Hurwitz--Kellogg discussion in Tabachnikov~\cite{tabachnikov2008}.
The interval weights needed here are derived directly below.

\subsection{Spectrum and Morse index rigidity}
The junction is allowed to move in the area variation problem, while
each outer edge stays on the fixed sphere. Competitors are not required
to remain minimal, conformal, or balanced. There are no volume
constraints and no junction energy. Section~\ref{sec:variations}
defines the normal form domain, Morse index, and nullity.
The junction compatibility follows the common velocity formulation
in Wang~\cite{wang2022curvature}. Tran~\cite{tran2020index} relates
the smooth free boundary index to the Jacobi problem with fixed
boundary and its Dirichlet-to-Neumann operator.

The index calculations in~\cite{matinpour2024morse} concern rotationally
symmetric $Y$-singular minimal surfaces with a circular junction. The low-index rigidity results
in~\cite{matinpour2025complete} concern a different, noncompact problem.

\begin{theorem}\label{thm:flat}
For the standard two-dimensional compact $Y$-cone in
$\overline{\B^n}$, with the variations and normal nullity of
Section~\ref{sec:variations},
\[
 \Ind(YC)=n-1,\qquad \nul(YC)=3(n-2).
\]
The coupled harmonic Steklov problem has eigenvalue $0$ with
multiplicity $n-1$. For every integer $k\geq1$, $k$ is an eigenvalue
of multiplicity $3(n-2)$. The eigenspace at $1$ is the normal kernel
of the area Hessian. Every kernel field is the normal component of an
ambient infinitesimal rotation.
\end{theorem}

\begin{theorem}\label{thm:rigidity}
For a disk-type free boundary minimal $Y$-surface $\Sigma$ satisfying
Definition~\ref{def:freeboundary},
\[
 \Ind(\Sigma)\geq n-1.
\]
Equality holds if and only if the three faces are the half-unit disks
of the standard flat $Y$-cone, each covered once, up to an orthogonal
transformation of $\R^n$. In particular, a nonflat disk-type
free boundary minimal $Y$-surface satisfying
Definition~\ref{def:freeboundary} has index at least $n$.
\end{theorem}

Its main identity is
\[
 Q(P^\perp v,P^\perp v)
 =-2\sum_{j=1}^3\int_{\Sigma_j}|P_j^\perp v|^2\,dA,
 \qquad v\in\R^n\text{ constant}.
\]
Fraser and Schoen~\cite[Theorem~3.1]{fraser2016sharp} prove the smooth
free boundary version. We check the junction flux and the endpoint
terms of the outer arcs. The rank argument giving the lower bound
and the stronger bound for nonflat surfaces does not use disk topology.
The theorem retains the common class of this paper, and its equality
case uses the planar classification in Lemma~\ref{lem:planar-triple}.
Section~\ref{sec:ambient-comparison} proves the same model counts and
index characterization for variations induced by ambient flows.

A related spectrum occurs in Wang's work on embedded closed stationary
triple-junction networks in $\Sph^2$~\cite{wang2024networks}.
Their length index is the number of complementary regions minus one,
and their nullity is three. For the theta network that is the link
of $YC\subset\R^3$, these numbers are two and three.
Theorem~\ref{thm:flat} instead concerns the area Hessian on the whole
surface, including its radial modes and additional normal directions.

\subsection{Scope and organization}
The three rigidity results use different additional hypotheses and do
not settle the unrestricted classification problem. Independently of
those hypotheses, Proposition~\ref{prop:area} gives
\[
 2\Area(\Sigma)=L
 =3\pi+\frac12\sum_j\int_{\Sigma_j}|\II_j|^2\,dA,
 \qquad L=\sum_j\Length\bigl(X_j(\sigma_j)\bigr),
\]
where each length is counted with its parametrization. This identity
characterizes equality in the area lower bound. It does not prove a
universal area upper bound.

Section~\ref{sec:setting} introduces the common setting and proves the
planar classification used by all three approaches.
Sections~\ref{sec:holomorphic}--\ref{sec:torsion} contain the conformal
and geometric proofs. Sections~\ref{sec:variations} and~\ref{sec:flat}
derive the Hessian and compute the flat spectrum.
Section~\ref{sec:translations} proves the translation identity and
index rigidity. Section~\ref{sec:area} gives the area identity and
the remaining question. The single appendix proves the conformal
boundary regularity used in the geometric arguments.

\subsection*{Acknowledgment}
I would like to thank Professor Jaigyoung Choe for suggesting the
problem and Professor Jacob Bernstein for helpful comments on an
earlier version of this manuscript.

\section{Geometric setting, regularity, and notation}\label{sec:setting}

\subsection{Abstract faces and their immersion}
Fix $0<\alpha<1$. Let $S_1,S_2,S_3$ be compact connected
$C^{2,\alpha}$ surfaces with corners. Each $S_j$ is a topological disk,
and its boundary consists of two regular $C^{2,\alpha}$ arcs,
\[
 \partial S_j=\gamma_j\cup\sigma_j.
\]
The interiors of the arcs are disjoint, and their endpoints are the
two marked vertices. Fix $C^{2,\alpha}$ diffeomorphisms
$\phi_j:I\to\gamma_j$, where $I=[-1,1]$. Gluing the three edges
by these maps gives the abstract triple junction.

Let $X_j:S_j\to\overline{\B^n}$ be $C^{2,\alpha}$ immersions with
rank two at every point, including the vertices. They have a common
junction map
\begin{equation}\label{eq:common-immersion}
 X_1\circ\phi_1=X_2\circ\phi_2=X_3\circ\phi_3
 =\gamma:I\to\overline{\B^n},\qquad |\gamma'|>0.
\end{equation}
At an edge or corner, $C^{2,\alpha}$ means restriction of a map of that
class from a neighborhood in an ordinary half-plane or quadrant chart.
The induced metrics are positive definite and $C^{1,\alpha}$, and
their coefficients extend locally across the coordinate boundaries.
These extensions are used only for the coordinate constructions below
and are not assumed to come from minimal immersions.

The face interiors and the open junction lie in $\B^n$, and
$X_j(\sigma_j)\subset\Sph^{n-1}$. We write $\Sigma_j$ for the
immersed face, $\Gamma$ for the common junction with its prescribed
parametrization, and $\Sigma$ for the triple. All integrals are taken
on the abstract faces with their induced metrics. This convention
counts multiplicity when images intersect.

\begin{definition}\label{def:freeboundary}
The triple is minimal if each face has zero mean curvature vector
and its outward unit junction conormals satisfy
\begin{equation}\label{eq:balance}
 \tau_1+\tau_2+\tau_3=0\qquad\text{on }\Gamma.
\end{equation}
It has spherical free boundary if the outward unit conormal $\eta_j$
to each outer edge equals the position vector,
\begin{equation}\label{eq:freeboundary}
 \eta_j=X_j\qquad\text{on }\sigma_j.
\end{equation}
A regular disk-type free boundary minimal $Y$-surface satisfies these
conditions and all the geometric assumptions above.
\end{definition}

We reserve $\tau_j$ for the junction conormal and $\eta_j$ for the
spherical boundary conormal. Since the $\tau_j$ are unit vectors,
\eqref{eq:balance} gives
\[
 \ip{\tau_i}{\tau_j}=-\frac12\qquad(i\ne j).
\]
Their span is two-dimensional, while the three face tangent planes
intersect along the line spanned by the unit tangent $T$ to the
junction.

The induced metrics $g_j=X_j^*g_{\R^n}$ agree on the abstract
junction after pullback,
\begin{equation}\label{eq:metric-pullback}
 \phi_j^*(g_j|_{T\gamma_j})=\gamma^*g_{\R^n}=:g_I.
\end{equation}
Define the geodesic curvature with the outward conormal by
\[
 \kappa_j=\ip{D_TT}{\tau_j}.
\]
Then
\begin{equation}\label{eq:kappa-balance}
 \sum_j\kappa_j=\ip{D_TT}{\sum_j\tau_j}=0.
\end{equation}

\subsection{The model and conformal boundary regularity}
We use the right half-disk
\begin{equation}\label{eq:half-disk}
 \begin{gathered}
 D_+=\{z=x+iy\in\C:x>0,\ |z|<1\},\\
 \gamma_0=\{iy:-1\leq y\leq1\},\qquad
 \sigma_0=\{z:|z|=1,\ \Rea z\geq0\}.
 \end{gathered}
\end{equation}
Here $\gamma_0$ is the model diameter and $\sigma_0$ the model semicircle.

\begin{definition}\label{def:Ycone}
The abstract model $YC$ consists of three copies of $\overline{D_+}$
with $\gamma_0$ identified pointwise. For its standard realization,
choose a unit vector $e$ and unit vectors $v_1,v_2,v_3\in e^\perp$
satisfying $v_1+v_2+v_3=0$. The face maps are
\[
 (x,y)\longmapsto xv_j+ye.
\]
Their union is the standard flat $Y$-cone in $\overline{\B^n}$.
Its linear span has dimension three.
\end{definition}

The next definition concerns regular conformal coordinates on each
face separately. It does not impose the matching condition of
Theorem~\ref{thm:main1}.

\begin{definition}\label{def:boundary-regular}
A surface as above satisfies the conformal boundary regularity
condition if each face admits a separate conformal parametrization
\[
 u_j=X_j\circ\psi_j:\overline{D_+}\longrightarrow\overline{\B^n}.
\]
The map $\psi_j$ is a homeomorphism onto $S_j$, is a diffeomorphism
on the interiors and open edges, and satisfies
\[
 \psi_j(\gamma_0)=\gamma_j,\qquad
 \psi_j(\sigma_0)=\sigma_j.
\]
We require $u_j\in C^{2,\alpha}(\overline{D_+})$ and
$\operatorname{rank}du_j=2$ everywhere on the closure.
\end{definition}

These symbols distinguish the abstract edges from the model edges
and their immersed images,
\[
 \gamma_0\xrightarrow{\ \psi_j\ }\gamma_j
       \xrightarrow{\ X_j\ }\Gamma,
 \qquad
 \sigma_0\xrightarrow{\ \psi_j\ }\sigma_j
       \xrightarrow{\ X_j\ }X_j(\sigma_j)\subset\Sph^{n-1}.
\]
Here the phrase model diameter refers to $\gamma_0$, rather than to
the general geometric junction.

Choose one orientation of the abstract junction and let
\[
 \gamma_\Gamma:[0,\ell]\longrightarrow\overline{\B^n}
\]
be its arclength parametrization. Each conformal face map determines
an increasing change of parameter
\[
 s_j:[-1,1]\longrightarrow[0,\ell].
\]
After orienting the face charts consistently, their traces satisfy
\begin{equation}\label{eq:separate-traces}
 \begin{gathered}
 u_j(iy)=\gamma_\Gamma(s_j(y)),\qquad
 \lambda_j(y)=s_j'(y)>0,\\
 s_j(-1)=0,\qquad s_j(1)=\ell.
 \end{gathered}
\end{equation}
The function $s_j$ converts the model coordinate $y$ to arclength on
the junction, and the three functions may differ. Under the matching
hypothesis of Theorem~\ref{thm:main1}, they agree through the
prescribed abstract gluing, although $y$ itself is not an arclength
parameter. The formulation on the abstract junction also applies when
its image has self-intersections.

Compactness and the nonvanishing differential give constants
\[
 0<\lambda_*\leq\lambda_j(y)\leq\lambda^*<\infty,
 \qquad -1\leq y\leq1.
\]
The second fundamental forms and the functions $c_j$ defined in
\eqref{eq:intro-c} are H\"older-continuous along the closed junction.

\begin{remark}\label{rem:regularity}
Proposition~\ref{prop:automatic-regularity} derives
Definition~\ref{def:boundary-regular} from the geometric regularity
in Definition~\ref{def:freeboundary} together with the right-angle
condition proved below. The resulting charts are separate and their
junction parameters are not identified; branch points and singular
vertices remain outside the class considered here.
\end{remark}

\begin{lemma}\label{lem:corners}
Let $\Sigma$ satisfy Definition~\ref{def:freeboundary}.
The junction meets the sphere radially at both endpoints, and each
face has an intrinsic angle $\pi/2$ there. Moreover,
\begin{equation}\label{eq:mixed-endpoint}
 \II_j(T,\tau_j)=0\qquad\text{at both endpoints}.
\end{equation}
\end{lemma}
\begin{proof}
At an endpoint with position vector $p\in\Sph^{n-1}$, the spherical
free boundary condition \eqref{eq:freeboundary} puts $p$ in every
face tangent plane. Their common intersection is $\R T$, so
\[
 p=\pm T.
\]
The outer boundary tangent is perpendicular to $p$, hence to $T$.
The two boundary rays of each face are orthogonal. Since the opposite
junction direction leaves the ball, the inward-pointing junction ray
is $-p$. The immersed quadrant lies in the supporting half-space of
the sphere. Its intrinsic angle is therefore $\pi/2$, not $3\pi/2$.

For a unit tangent $Z$ to the outer edge, differentiation of
$\eta_j=X_j$ gives
\begin{equation}\label{eq:mixed-fb}
 \II_j(Z,\eta_j)=(D_Z\eta_j)^\perp=(D_ZX_j)^\perp=0.
\end{equation}
At a vertex, $\eta_j=\pm T$ and $Z$ is parallel to $\tau_j$.
The continuous endpoint limit proves \eqref{eq:mixed-endpoint}.
\end{proof}

\subsection{Counting sign changes}
\begin{definition}\label{def:signchanges}
For a continuous real function $f$ on the oriented junction,
$\nsc(f)$ is the supremum of the integers $k\geq0$ for which there
are points in the open junction satisfying
\[
 s_0<\cdots<s_k,\qquad f(s_i)\ne0,\qquad
 f(s_i)f(s_{i+1})<0\quad(0\leq i<k).
\]
Set $\nsc(0)=0$, and allow $\nsc(f)=+\infty$.
Zeros without a sign change, including intervals of zeros, do not
increase this count. It is unchanged by reversal of orientation,
increasing reparametrization, or multiplication by a continuous
factor that is nowhere zero.
\end{definition}

\subsection{Planar classification}
The next lemma identifies the surface once all three faces are
planar. Lemma~\ref{lem:one-planar} later shows, using the holomorphic
boundary identities, that one planar face already forces the other
two to be planar.

\begin{lemma}\label{lem:planar-triple}
Let $\Sigma$ satisfy Definition~\ref{def:freeboundary}.
If all three faces are planar, each image $X_j(S_j)$ is a half-unit
disk covered once. Together the three images form the standard
flat $Y$-cone, up to an orthogonal transformation of $\R^n$.
\end{lemma}
\begin{proof}
\proofstep{1}{Identify the face planes and the junction.}
Suppose a face lies in the affine plane $a_j+P_j$, where $P_j$ is a
linear two-plane. At an outer boundary point $p$, free boundary gives
\[
 p=\eta_j\in P_j,\qquad p\in a_j+P_j.
\]
Since $p$ belongs to both $P_j$ and $a_j+P_j$, we have $a_j\in P_j$.
The affine plane is therefore the linear plane $P_j$.
Two face planes cannot coincide, since their junction conormals
would be collinear, contrary to
$\ip{\tau_i}{\tau_j}=-1/2$.

The common regular junction lies in their intersection, a line $L$
through the origin. Choose a unit vector $e\in L$ and write
\[
 \gamma(t)=r(t)e.
\]
Since $r'$ never vanishes, $r$ is strictly monotone. Its endpoint
images lie on the sphere and must therefore be the two opposite ends
of $L\cap\overline{\B^n}$, which shows that the full diameter is
covered once.

\proofstep{2}{Identify each face and its multiplicity.}
Double $S_j$ along $\gamma_j$ and reflect the image of the second
copy in $L$ within $P_j$. This gives the continuous map
\[
 F_j:\widehat S_j\longrightarrow\overline{D_j},
 \qquad D_j=P_j\cap\B^n,
\]
where $\widehat S_j$ is a topological disk and its boundary maps to
$\partial D_j$. Away from the doubled edge, $F_j$ is a local
homeomorphism because $X_j$ is a planar immersion. At an interior
point of that edge, the immersion gives a half-plane coordinate.
Since the face lies in the fixed plane $P_j$, the unit conormal
$\tau_j$ is perpendicular to $L$ in $P_j$. Along the connected
junction it is one of the two fixed unit normals to $L$ in $P_j$.
Continuity fixes that choice, so the face stays locally on the same
side of $L$. The reflected copy lies on the other side, and together
the two copies map homeomorphically onto a full neighborhood.
The two edge endpoints belong to $\partial\widehat S_j$.

The restriction to the interior is proper. Indeed, a sequence whose
images stay in a compact subset of $D_j$ cannot approach
$\partial\widehat S_j$, whose image lies on $\partial D_j$.
Its image is open by local invertibility and closed by properness,
so it is all of $D_j$. A proper local homeomorphism is a covering.
The interior of $\widehat S_j$ is connected and $D_j$ is simply connected,
so this covering has one sheet.

The doubled edge already maps once onto the whole diameter.
Injectivity excludes any other preimage of that diameter. The two
components of its complement map to the two half-disks in $D_j$.
Each original face is therefore a full half-unit disk covered once.

\proofstep{3}{Use the conormal balance.}
The transverse unit directions of the three half-disks sum to zero.
They span a two-plane perpendicular to $L$ and meet at $120^\circ$,
which gives the standard realization of $YC$.
\end{proof}

\section{Holomorphic boundary identities}\label{sec:holomorphic}
We extend the Euclidean inner product complex bilinearly to $\C^n$;
it is not the Hermitian product. For a conformal minimal face map
$u$ satisfying Definition~\ref{def:boundary-regular}, write
\[
 u_z=\frac12(u_x-iu_y),\qquad
 \Hopf=u_{zz}^{\perp},\qquad
 \varphi=\ip{\Hopf}{\Hopf}.
\]
The coefficient $\Hopf$ takes values in the complexified normal bundle,
whereas $\varphi$ is a complex scalar. The associated differentials
are $\Hopf\,dz^2$ and $\varphi\,dz^4$.

\begin{lemma}\label{lem:quartic}
The function $\varphi$ is holomorphic on $D_+$ and belongs to
$C^{0,\alpha}(\overline{D_+})$. On the model semicircle $\sigma_0$,
\begin{equation}\label{eq:quartic-positive}
 z^4\varphi(z)=\frac14|u_{rr}^{\perp}|^2\geq0.
\end{equation}
\end{lemma}
\begin{proof}
Harmonicity and conformality give
\[
 u_{z\bar z}=0,\qquad \ip{u_z}{u_z}=0.
\]
It follows that $u_{zz}$ is a $\C^n$-valued holomorphic function and
that
\[
 \ip{u_{zz}}{u_z}=0.
\]
The complex tangent plane is spanned by $u_z,u_{\bar z}$, with
$\ip{u_z}{u_{\bar z}}>0$. Its tangential projection formula is
\[
 u_{zz}=u_{zz}^{\perp}
 +\frac{\ip{u_{zz}}{u_{\bar z}}}{\ip{u_z}{u_{\bar z}}}\,u_z.
\]
The normal term is orthogonal to $u_z$, and $\ip{u_z}{u_z}=0$.
Squaring gives
\begin{equation}\label{eq:square-identity}
 \varphi=\ip{u_{zz}^{\perp}}{u_{zz}^{\perp}}
        =\ip{u_{zz}}{u_{zz}}.
\end{equation}
The last expression is holomorphic, and hence so is $\varphi$. Only
the scalar square is used here; no holomorphicity is claimed for
$u_{zz}^{\perp}$ as a map into a fixed copy of $\C^n$. This quartic
construction also appears in Fraser and
Schoen~\cite{fraser2015uniqueness}.
The regularity of $u$ and its nonvanishing differential give the
stated H\"older extension of the projections and of $\varphi$.

On the spherical boundary $r=1$, conformality makes $u_r$ point in
the outward conormal direction. Free boundary therefore gives
\[
 u_r=\rho u,\qquad \rho=|u_r|>0.
\]
Differentiating along the boundary,
\[
 u_{r\theta}=\rho_\theta u+\rho u_\theta.
\]
Both $u$ and $u_\theta$ are tangent to the face there, so
\[
 u_{r\theta}^{\perp}=0\qquad\text{on }\sigma_0.
\]
For $r>0$, the polar expression for the complex derivative is
\[
 u_{zz}=\frac{e^{-2i\theta}}4
 \left(u_{rr}-\frac1r u_r-\frac{2i}{r}u_{r\theta}
       -\frac1{r^2}u_{\theta\theta}+\frac{2i}{r^2}u_\theta\right).
\]
After normal projection, the terms involving $u_r$ and $u_\theta$
disappear, and we obtain
\[
 \Hopf=\frac{e^{-2i\theta}}4
 \left(u_{rr}^{\perp}-\frac1{r^2}u_{\theta\theta}^{\perp}
       -\frac{2i}{r}u_{r\theta}^{\perp}\right).
\]
Harmonicity in polar coordinates reads
\[
 u_{rr}+\frac1r u_r+\frac1{r^2}u_{\theta\theta}=0.
\]
Taking the normal part at $r=1$ gives
$u_{rr}^{\perp}+u_{\theta\theta}^{\perp}=0$, and therefore
\[
 z^2\Hopf=\frac12u_{rr}^{\perp}\qquad\text{on }\sigma_0.
\]
The right side is a real normal vector, so its complex-bilinear
square is an ordinary nonnegative norm and gives
\eqref{eq:quartic-positive}. The identity is first obtained on the
open semicircle and then extended to the endpoints by $C^{2,\alpha}$
regularity.
\end{proof}

\begin{lemma}\label{lem:planarity}
For a face map as above, the following statements hold.
\begin{enumerate}[label=\textup{(\roman*)},leftmargin=2.4em]
\item If its second fundamental form vanishes along the open
spherical boundary arc, the entire face lies in a linear two-plane
through the origin.
\item If $\Ima\varphi=0$ on the model diameter $\gamma_0$, the
second fundamental form vanishes along the open spherical boundary
arc. In particular, the face is planar.
\end{enumerate}
\end{lemma}
\begin{proof}
For (i), parametrize the outer arc by arclength as $c(s)$.
The tangent plane of the face is spanned by $c$ and $c'$, since its
outward conormal is $c$. Vanishing of the second fundamental form
makes $c''$ tangent. The identities
\[
 |c|=|c'|=1,\qquad
 \ip{c''}{c}=-1,\qquad \ip{c''}{c'}=0
\]
force $c''=-c$, so the boundary is a great-circle arc contained in a
linear two-plane $P$.

For every constant $a\in P^\perp$, the function
\[
 w=\ip{a}{u}
\]
is harmonic. On the open outer arc, free boundary gives
\[
 w=0,\qquad \partial_rw=\rho w=0.
\]
To see that these data imply $w=0$, straighten a short circular
subarc conformally and reflect the zero-Dirichlet harmonic function
oddly across the resulting interval. Its complex derivative vanishes
on the interval, since both tangential and normal derivatives are
zero. The identity theorem gives local vanishing, and harmonic unique
continuation extends it throughout $D_+$. Applying this to a basis
of $P^\perp$ puts the entire face in $P$.

For (ii), set $H=z^4\varphi$. On $\gamma_0$ we have $z^4=y^4\in\R$,
so the hypothesis makes $H$ real on the diameter;
Lemma~\ref{lem:quartic} gives the same conclusion on $\sigma_0$.
The maximum principle gives $\Ima H=0$ on $D_+$. Since a real-valued
holomorphic function is constant and $H(0)=0$ by the boundedness of
$\varphi$ at the boundary point $0$, we have $H\equiv0$.

By \eqref{eq:quartic-positive}, $u_{rr}^{\perp}=0$ on the open
semicircle. The preceding boundary calculation also gives
\[
 u_{r\theta}^{\perp}=0,\qquad
 u_{\theta\theta}^{\perp}=-u_{rr}^{\perp}=0.
\]
Because $u_r$ and $u_\theta$ span the tangent plane, these are all
components of the second fundamental form and part (i) applies. The
real norm in \eqref{eq:quartic-positive} is important here: in higher
codimension, vanishing of a complex-bilinear square by itself would
not force the underlying vector to vanish.
\end{proof}

\section{Proof of conformal rigidity}\label{sec:conformal-proof}
\begin{proof}[Proof of Theorem~\ref{thm:main1}]
\proofstep{1}{Use the common junction parameter.}
First the face interiors lie in $\B^n$, since
\[
 \Delta|u_j|^2=2\bigl(|(u_j)_x|^2+|(u_j)_y|^2\bigr)>0.
\]
An open junction point cannot map to the sphere. If its common image
were $p\in\Sph^{n-1}$, the Hopf boundary lemma for $|u_j|^2$ would
give a strictly positive outward derivative on each model diameter.
Conformality and $\nu=-\partial_x$ give
\[
 0<\partial_\nu|u_j|^2
 =2|(u_j)_y|\ip{p}{\tau_j}.
\]
All three inner products would be positive, contradicting
$\sum_j\tau_j=0$. The maps therefore satisfy the containment
conditions of Section~\ref{sec:setting}.

Write their common trace as
\[
 v(y)=u_1(0,y)=u_2(0,y)=u_3(0,y).
\]
Along the model diameter,
\[
 (u_j)_y=v',\qquad (u_j)_{yy}=v'',\qquad
 |(u_j)_x|=|v'|>0.
\]
Since $\tau_j=-(u_j)_x/|v'|$, the balance condition and its derivative
with respect to the common parameter give
\begin{equation}\label{eq:cartesian-balance}
 \sum_j(u_j)_x=0,\qquad \sum_j(u_j)_{xy}=0
 \qquad\text{on }\gamma_0.
\end{equation}
Here the matching hypothesis enters: the three faces share the same
parameter $y$ along the diameter, so the differentiated balance in
\eqref{eq:cartesian-balance} is meaningful.

\proofstep{2}{Apply the maximum principle to the summed quartic coefficient.}
By Lemma~\ref{lem:quartic}, the functions
\[
 h(z)=\sum_j\varphi_j(z)
     =\sum_j\ip{(u_j)_{zz}}{(u_j)_{zz}},\qquad H(z)=z^4h(z)
\]
are holomorphic in $D_+$ and continuous on its closure.
On the diameter, harmonicity gives $(u_j)_{xx}=-v''$. Hence
\[
 (u_j)_{zz}=-\frac12\bigl(v''+i(u_j)_{xy}\bigr),
\]
and the complex-bilinear square is
\[
 \varphi_j=\frac14\left(|v''|^2-|(u_j)_{xy}|^2
                   +2i\ip{v''}{(u_j)_{xy}}\right).
\]
Equation~\eqref{eq:cartesian-balance} gives
\[
 \Ima h=\frac12\ip{v''}{\sum_j(u_j)_{xy}}=0
 \qquad\text{on }\gamma_0.
\]
Since $(iy)^4$ is real, $H$ is real on the diameter.
On the outer semicircle, Lemma~\ref{lem:quartic} gives
\[
 H=\frac14\sum_j|(u_j)_{rr}^{\perp}|^2\geq0.
\]
The imaginary part of $H$ therefore has zero values on the entire
boundary. The maximum principle makes $H$ real throughout $D_+$,
so it is constant. Continuity at $0$ gives $H(0)=0$, and hence
$H\equiv0$.

\proofstep{3}{Conclude planarity and identify the image.}
On the semicircle, the sum of real squares above can vanish only if
each term vanishes. For every face,
\[
 (u_j)_{rr}^{\perp}=0,\qquad
 (u_j)_{r\theta}^{\perp}=0,\qquad
 (u_j)_{\theta\theta}^{\perp}=0.
\]
The last two identities use free boundary and normal harmonicity
as in Lemma~\ref{lem:quartic}. Since $(u_j)_r,(u_j)_\theta$ span
the tangent plane, $\II_j=0$ along its open spherical boundary arc.
Lemma~\ref{lem:planarity}(i) makes every face planar.
Lemma~\ref{lem:planar-triple} identifies the full half-unit disks,
each covered once, and gives the standard $Y$-cone.
\end{proof}

\section{Rigidity without matching conformal parametrizations}\label{sec:oscillation}
Let $\Sigma$ satisfy Definition~\ref{def:freeboundary}. The right
angles from Lemma~\ref{lem:corners}, together with
Proposition~\ref{prop:automatic-regularity}, give separate conformal
maps $u_j$ satisfying Definition~\ref{def:boundary-regular}. Their
junction traces have the form \eqref{eq:separate-traces}, with no
identification among the three functions $s_j$.

\subsection{Geometric junction data}
\begin{lemma}\label{lem:c-balance}
Along the geometric junction, set
\[
 A_j=\II_j(T,T),\qquad B_j=\II_j(T,\tau_j),\qquad
 c_j=\ip{A_j}{B_j}.
\]
Then the pointwise identity
\begin{equation}\label{eq:c-balance}
 c_1+c_2+c_3=0
\end{equation}
holds on $\Gamma$.
For each separate face map, let
\[
 \Hopf_j=(u_j)_{zz}^{\perp},\qquad
 \varphi_j=\ip{\Hopf_j}{\Hopf_j},\qquad
 g_j(y)=\Ima\varphi_j(iy).
\]
With $s_j$ and $\lambda_j=s_j'>0$ from \eqref{eq:separate-traces},
\begin{equation}\label{eq:quartic-trace}
 \begin{aligned}
 \Hopf_j(iy)
 &=-\frac{\lambda_j(y)^2}{2}
       \bigl(A_j(s_j(y))-iB_j(s_j(y))\bigr),\\
 \varphi_j(iy)
 &=\frac{\lambda_j(y)^4}{4}
       \bigl(|A_j(s_j(y))|^2-|B_j(s_j(y))|^2
                            -2i c_j(s_j(y))\bigr),\\
 g_j(y)&=-\frac{\lambda_j(y)^4}{2}\,c_j(s_j(y)).
 \end{aligned}
\end{equation}
Here $A_j(s)$ denotes the geometric vector at $\gamma_\Gamma(s)$,
and similarly for $B_j$ and $c_j$.
\end{lemma}
\begin{proof}
Write $K=D_TT$ and $\kappa_j=\ip{K}{\tau_j}$.
The Gauss formula in the orthonormal tangent frame $(T,\tau_j)$ gives
\[
 A_j=K-\kappa_j\tau_j,\qquad
 B_j=D_T\tau_j+\kappa_jT.
\]
The relations $K\perp T$, $\tau_j\perp T$, and
$\ip{\tau_j}{D_T\tau_j}=0$ imply
\[
 \begin{aligned}
 c_j
 &=\ip{K-\kappa_j\tau_j}{D_T\tau_j+\kappa_jT}\\
 &=\ip{K}{D_T\tau_j}.
 \end{aligned}
\]
Differentiating $\sum_j\tau_j=0$ along the common tangent gives
\[
 \sum_jc_j=\ip{K}{\sum_jD_T\tau_j}=0.
\]

We now express the same quantities in the half-disk coordinates.
From \eqref{eq:separate-traces} and arclength parametrization,
\[
 (u_j)_y(iy)=\lambda_j(y)T.
\]
Conformality gives $|(u_j)_x|=\lambda_j$ and orthogonality of the
coordinate vectors. Since $-\partial_x$ is outward along the model
diameter,
\[
 (u_j)_x(iy)=-\lambda_j(y)\tau_j.
\]
Normal projection of the second derivatives therefore gives
\[
 \bigl((u_j)_{yy}\bigr)^\perp=\lambda_j^2 A_j,\qquad
 \bigl((u_j)_{xx}\bigr)^\perp=-\lambda_j^2 A_j,\qquad
 \bigl((u_j)_{xy}\bigr)^\perp=-\lambda_j^2 B_j,
\]
where the right sides are evaluated at $s_j(y)$. The second identity
uses minimality, namely $\II_j(\tau_j,\tau_j)=-A_j$.
Substitution into
\[
 (u_j)_{zz}=\frac14\bigl((u_j)_{xx}-(u_j)_{yy}-2i(u_j)_{xy}\bigr)
\]
proves the first line of \eqref{eq:quartic-trace}. Taking the
complex-bilinear square gives the second line of
\eqref{eq:quartic-trace}, and its imaginary part gives the third.
The formulas extend to the endpoints by continuity.
\end{proof}

\subsection{A scalar endpoint formula}
The following lemma expresses the value at a corner in terms of the
imaginary boundary values on the model diameter.

\begin{lemma}\label{lem:endpoint}
Let $m\in\{1,2\}$, and let $F$ be holomorphic on $D_+$ and
$C^{0,\alpha}$ on its closure. Suppose $z^{2m}F(z)$ is real on
$\sigma_0$, and set $g(y)=\Ima F(iy)$. Then $g(\pm1)=0$,
$F(i)$ is real, and
\begin{equation}\label{eq:endpoint-general}
 F(i)=\frac1\pi\int_{-1}^{1}
             \frac{1+y^{4m-1}}{1-y}\,g(y)\,dy.
\end{equation}
The integral is absolutely convergent.
\end{lemma}
\begin{proof}
\proofstep{1}{Reflect across the semicircle.}
The function $H=z^{2m}F$ is real on the semicircle. Reflecting it
across the unit circle gives $\widetilde H(z)=\overline{H(1/\bar z)}$,
and division by $z^{2m}$ produces the following holomorphic extension
to the right half-plane:
\[
 \widetilde F(z)=
 \begin{cases}
 F(z),&|z|\leq1,\\
 z^{-4m}\overline{F(1/\bar z)},&|z|\geq1.
 \end{cases}
\]
The formulas agree on the open semicircle, so Schwarz reflection
makes $\widetilde F$ holomorphic on the whole open right half-plane.
Since $\pm i$ are boundary points of the right half-plane, the
continuation introduces no interior singularity there. The extension is
bounded, continuous to the imaginary axis, and $O(|z|^{-4m})$ at
infinity. Since $(\pm i)^{2m}=(-1)^m$, continuity of the original
boundary condition also gives $F(\pm i)\in\R$ and $g(\pm1)=0$.

The imaginary boundary data of the extension are
\[
 G(t)=\Ima\widetilde F(it)=
 \begin{cases}
 g(t),&|t|\leq1,\\
 -t^{-4m}g(1/t),&|t|\geq1.
 \end{cases}
\]
In particular, $G$ is bounded and integrable. It is
H\"older-continuous near $\pm1$ and vanishes there.

\proofstep{2}{Use the half-plane Poisson formula.}
For $\Rea z>0$, define
\[
 S(z)=\frac{i}{\pi}\int_\R\frac{G(t)}{z-it}\,dt.
\]
This integral defines a holomorphic function. Its imaginary part is
\[
 \Ima S(x+iy)=\frac1\pi\int_\R
       \frac{x}{x^2+(y-t)^2}G(t)\,dt,
\]
the Poisson extension of $G$. Uniqueness for bounded harmonic
functions on the half-plane gives
$\Ima S=\Ima\widetilde F$. Their difference is therefore a real
constant. Along the positive real axis,
\[
 |S(x)|\leq\frac{\|G\|_{L^1}}{\pi x}\longrightarrow0,
 \qquad \widetilde F(x)=O(x^{-4m}),
\]
so that constant is zero.

\proofstep{3}{Take the endpoint limit.}
For $x>0$,
\[
 \Rea\widetilde F(x+i)=\frac1\pi\int_\R
       \frac{1-t}{x^2+(1-t)^2}G(t)\,dt.
\]
Near $t=1$, H\"older continuity gives
\[
 \left|\frac{1-t}{x^2+(1-t)^2}G(t)\right|
 \leq\frac{|G(t)|}{|1-t|}
 \leq C|1-t|^{\alpha-1}.
\]
Since $\alpha>0$, this singularity is integrable. Away from $1$, the
boundedness and decay of $G$ provide the required majorant, so
dominated convergence yields
\[
 F(i)=\frac1\pi\int_\R\frac{G(t)}{1-t}\,dt.
\]
On each of the two intervals $|t|>1$, substitute $t=1/y$.
Using the reflected data and reversing the limits gives
\[
 \int_{|t|>1}\frac{G(t)}{1-t}\,dt
 =\int_{-1}^{1}\frac{y^{4m-1}}{1-y}g(y)\,dy.
\]
Adding the contribution from $[-1,1]$ proves
\eqref{eq:endpoint-general}. The same endpoint estimate proves its
absolute convergence.
\end{proof}

\subsection{Three quartic moment identities}
\begin{lemma}\label{lem:moments}
Suppose $F$ satisfies Lemma~\ref{lem:endpoint} with $m=2$, and let
$g(y)=\Ima F(iy)$. Set
\begin{equation}\label{eq:Wt}
 W(y)=(1-y^2)(1+y^2)^2,\qquad t(y)=\frac{y}{1+y^2}.
\end{equation}
Then $W>0$ on $(-1,1)$, and $t$ is strictly increasing from
$[-1,1]$ onto $[-1/2,1/2]$. Moreover,
\begin{equation}\label{eq:three-moments}
 \int_{-1}^{1}W(y)t(y)^k g(y)\,dy=0,\qquad k=0,1,2.
\end{equation}
\end{lemma}
\begin{proof}
Choose
\[
 P_0(z)=1+z^6,\qquad P_1(z)=i(z-z^5),\qquad P_2(z)=z^2+z^4.
\]
These polynomials make the following contour integrands real on the
semicircle. Indeed, if $z=e^{i\theta}$, then
\[
 iF(z)P_k(z)\,dz=-(z^4F(z))\frac{P_k(z)}{z^3}\,d\theta,
\]
and
\[
 \frac{P_0}{z^3}=2\cos3\theta,\qquad
 \frac{P_1}{z^3}=2\sin2\theta,\qquad
 \frac{P_2}{z^3}=2\cos\theta.
\]
Since $z^4F$ is real there, each semicircle integral has zero
imaginary part.

Cauchy's theorem gives
\[
 \int_{\partial D_+}iF(z)P_k(z)\,dz=0.
\]
For a coefficient continuous to the piecewise smooth boundary,
this follows by approximation with interior contours.
The positive boundary orientation traverses the diameter from $i$
to $-i$. Setting $z=iy$ and reversing the limits gives
\[
 \int_{\gamma_0}iF(z)P_k(z)\,dz
 =\int_{-1}^{1}F(iy)P_k(iy)\,dy.
\]
The diameter values are
\[
 P_0(iy)=1-y^6,\qquad P_1(iy)=-(y-y^5),\qquad
 P_2(iy)=-(y^2-y^4).
\]
Taking imaginary parts and changing the signs of the last two
zero identities, we obtain
\begin{equation}\label{eq:raw-moments}
 \begin{aligned}
 \int_{-1}^{1}(1-y^6)g(y)\,dy&=0,\\
 \int_{-1}^{1}(y-y^5)g(y)\,dy&=0,\\
 \int_{-1}^{1}(y^2-y^4)g(y)\,dy&=0.
 \end{aligned}
\end{equation}
Now
\[
 W=(1-y^6)+(y^2-y^4),\qquad Wt=y-y^5,\qquad Wt^2=y^2-y^4,
\]
which proves \eqref{eq:three-moments}. Finally,
\[
 t'(y)=\frac{1-y^2}{(1+y^2)^2}>0\quad(-1<y<1),\qquad
 t(\pm1)=\pm\frac12.
\]
\end{proof}

\Needspace{8\baselineskip}
\begin{lemma}\label{lem:analytic-oscillation}
Let $F$ satisfy Lemma~\ref{lem:endpoint} with $m=2$, and suppose
$z^4F(z)\geq0$ on $\sigma_0$. Set $g(y)=\Ima F(iy)$.
\begin{enumerate}[label=\textup{(\roman*)},leftmargin=2.4em]
\item If $\nsc(g)\leq2$, then $g\equiv0$ and $F\equiv0$.
\item If $g\not\equiv0$ and $\nsc(g)=3$, its successive nonzero
signs, read from $-1$ to $1$, have the pattern $-,+,-,+$.
\end{enumerate}
Zeros and intervals of zeros are ignored when reading this pattern.
\end{lemma}
\begin{proof}
The moment identities say that
\begin{equation}\label{eq:polynomial-moments}
 \int_{-1}^{1}W(y)p(t(y))g(y)\,dy=0
 \qquad\text{whenever }\deg p\leq2.
\end{equation}
Suppose $g\not\equiv0$ and $q=\nsc(g)\leq2$.
Choose a separating zero $\xi_r$ between each pair of successive
opposite-sign regions, with $\xi_1<\cdots<\xi_q$.
Such choices are possible because there are only finitely many
sign alternations. If a separating set is a zero interval, choose
any point in it. Since $t$ is increasing, a suitable overall sign
of the polynomial
\[
 p(t)=\prod_{r=1}^{q}(t-t(\xi_r))
\]
makes $p(t(y))g(y)$ nonnegative everywhere and positive on a nonempty
open interval. For $q=0$, use a constant of the appropriate sign.
Positivity of $W$ contradicts \eqref{eq:polynomial-moments}.
Thus $g=0$.

The function $z^4F$ is then real on both boundary arcs. The maximum
principle and holomorphicity make it constant. Its value at the
boundary point $0$ is zero, so $F\equiv0$.

When there are exactly three sign changes, choose separating values
$t_1<t_2<t_3$ in $(-1/2,1/2)$.
Lemma~\ref{lem:endpoint} can be rewritten as
\begin{equation}\label{eq:endpoint-R}
 F(i)=\frac1\pi\int_{-1}^{1}W(y)R(t(y))g(y)\,dy,
\end{equation}
where
\[
 R(t)=-\frac{t^2}{2}+\frac{3t}{4}+\frac78+\frac1{8(1-2t)}.
\]
Indeed, $W(y)R(t(y))=(1+y^7)/(1-y)$.
Let $p_2$ interpolate $R$ at $t_1,t_2,t_3$. Equation
\eqref{eq:polynomial-moments} allows us to subtract $p_2$ in
\eqref{eq:endpoint-R}.

The interpolation remainder has the exact form
\begin{equation}\label{eq:exact-remainder}
 R(t)-p_2(t)=
 \frac{(t-t_1)(t-t_2)(t-t_3)}
      {(1-2t)\prod_{r=1}^{3}(1-2t_r)}.
\end{equation}
To verify it, multiply the left side by $1-2t$.
The result extends to a polynomial of degree at most three that
vanishes at $t_1,t_2,t_3$. Its value at $1/2$ is $1/8$, which
determines the displayed factor. All denominator factors are positive
in the interval under consideration. Hence the remainder has
successive signs $-,+,-,+$.

If $g$ had the opposite pattern $+,-,+,-$, the integrand in
\[
 F(i)=\frac1\pi\int_{-1}^{1}
 W(y)\bigl(R(t(y))-p_2(t(y))\bigr)g(y)\,dy
\]
would be nonpositive everywhere and negative on a nonempty open set.
The integral is absolutely convergent by Lemma~\ref{lem:endpoint},
and the subtracted polynomial term is continuous. We would obtain
$F(i)<0$. This contradicts $F(i)=i^4F(i)\geq0$, the endpoint limit
of the assumed boundary inequality. Only $-,+,-,+$ is possible.
\end{proof}

\subsection{Completion of the geometric theorem}
\begin{lemma}\label{lem:one-planar}
Let $\Sigma$ satisfy Definition~\ref{def:freeboundary}.
If one face is planar, the other two faces are planar as well.
\end{lemma}
\begin{proof}
Suppose face $1$ is planar. Then $A_1=B_1=0$ on $\Gamma$, and the
decompositions in Lemma~\ref{lem:c-balance} give
\[
 D_TT=\kappa_1\tau_1,\qquad D_T\tau_1=-\kappa_1T.
\]
For $j=2,3$, differentiate the constant-angle identity to obtain
\[
 0=D_T\ip{\tau_1}{\tau_j}
   =\ip{D_T\tau_1}{\tau_j}+\ip{\tau_1}{D_T\tau_j}.
\]
The first term is $-\kappa_1\ip{T}{\tau_j}=0$. Therefore
\[
 c_j=\ip{D_TT}{D_T\tau_j}
     =\kappa_1\ip{\tau_1}{D_T\tau_j}=0.
\]
For each of the two separate face maps, \eqref{eq:quartic-trace} gives
\[
 \Ima\varphi_j(iy)
 =-\frac{\lambda_j(y)^4}{2}c_j(s_j(y))=0.
\]
Lemma~\ref{lem:planarity}(ii) makes both faces planar.
\end{proof}

\begin{proof}[Proof of Theorem~\ref{thm:oscillation}]
Equation~\eqref{eq:quartic-trace}, the positivity of $\lambda_j$,
and the increasing parametrization $s_j$ give
\[
 \nsc(g_j)=\nsc(c_j).
\]
For (i), Lemma~\ref{lem:analytic-oscillation}(i) makes $g_j=0$ on
one face. That face is planar by Lemma~\ref{lem:planarity}.
Lemma~\ref{lem:one-planar} makes the other faces planar, and
Lemma~\ref{lem:planar-triple} gives $YC$.

For (ii), suppose the surface is not $YC$. Part (i) then forces all
three counts to be exactly three, and
Lemma~\ref{lem:analytic-oscillation}(ii) gives the pattern $-,+,-,+$
for each $g_j$. The minus sign in
\eqref{eq:quartic-trace} gives $+,-,+,-$ for every $c_j$ in the
same arclength direction on $\Gamma$.

For each $j$, let
\[
 a_j=\inf\{s\in[0,\ell]:c_j(s)<0\}.
\]
The stated sign pattern implies $a_j>0$. Before $a_j$, the function
$c_j$ is nonnegative and is positive somewhere.
Choose $j_*$ for which $a_{j_*}=\min_j a_j$, and choose
$r<a_{j_*}$ with $c_{j_*}(r)>0$. All three values $c_j(r)$ are
nonnegative, and one is positive. This contradicts
$\sum_jc_j(r)=0$ from Lemma~\ref{lem:c-balance}.
The argument allows intervals of zeros.
\end{proof}

\begin{proof}[Proof of Corollary~\ref{cor:geometric}]
One planar face gives the conclusion by Lemmas~\ref{lem:one-planar}
and~\ref{lem:planar-triple}. If $B_j=0$ on the junction of one face,
then $c_j=\ip{A_j}{B_j}=0$. This and the stated alternative $c_j=0$
are covered by Theorem~\ref{thm:oscillation}(i).
For a nonflat surface, (i) forces at least three sign changes in
every $c_j$. If all counts were at most three, (ii) would make the
surface flat. At least one count is therefore four or larger.
\end{proof}

\begin{remark}
The moment and endpoint arguments are facewise and use the separate
parameters $s_j$ and factors $\lambda_j$. The faces are coupled only
through the junction identities, in particular $\sum_j c_j=0$, and
through Lemma~\ref{lem:one-planar} once one face is planar.
\end{remark}

\section{The three-dimensional torsion criterion}\label{sec:torsion}
Let $\Sigma\subset\overline{\B^3}$ satisfy
Definition~\ref{def:freeboundary}. By Lemma~\ref{lem:corners} and
Proposition~\ref{prop:automatic-regularity}, we may choose separate
conformal face maps with
\[
 u_j(iy)=\gamma_\Gamma(s_j(y)),
 \qquad \lambda_j(y)=s_j'(y)>0.
\]
The functions $s_j$ may differ from one face to another. In this
section we work with the scalar Hopf coefficients rather than their
quartic squares.

Along the junction, choose
\[
 N_j=T\times\tau_j,
 \qquad
 a_j=\ip{\II_j(T,T)}{N_j},
 \qquad
 b_j=\ip{\II_j(T,\tau_j)}{N_j}.
\]
The conormal balance implies $\sum_jN_j=0$. The Weingarten equation is
\[
 D_TN_j=-a_jT-b_j\tau_j.
\]
Differentiating the sum and separating its components parallel and
perpendicular to $T$ gives
\[
 \sum_ja_j=0,
 \qquad \sum_jb_j\tau_j=0.
\]
Since $\tau_3=-\tau_1-\tau_2$ and $\tau_1,\tau_2$ are linearly
independent, the second equality becomes
\[
 (b_1-b_3)\tau_1+(b_2-b_3)\tau_2=0,
\]
which forces
\begin{equation}\label{eq:ab-balance}
 a_1+a_2+a_3=0,
 \qquad b_1=b_2=b_3=:b.
\end{equation}
This is also the constant-angle geodesic-torsion relation discussed
by Brander and L\'opez~\cite[Lemma~2]{brander2017remarks}.
The equalities hold on the closed junction by continuity.

The unit normal on each face is
\[
 N_j(z)=\frac{(u_j)_x\times(u_j)_y}
 {|(u_j)_x\times(u_j)_y|}.
\]
On the model diameter, $(u_j)_x=-\lambda_j\tau_j$ and
$(u_j)_y=\lambda_jT$, so this normal agrees with $T\times\tau_j$.
Define
\[
 q_j(z)=\ip{(u_j)_{zz}}{N_j}.
\]
To verify that $q_j$ is holomorphic, harmonicity gives
\[
 \partial_{\bar z}q_j
 =\ip{(u_j)_{zz}}{(N_j)_{\bar z}}.
\]
The vector $(N_j)_{\bar z}$ is tangent to the complexified face.
Differentiating $\ip{N_j}{(u_j)_z}=0$ and using
$(u_j)_{z\bar z}=0$ gives
\[
 \ip{(N_j)_{\bar z}}{(u_j)_z}=0.
\]
In the complex tangent plane this forces $(N_j)_{\bar z}$ to be
proportional to $(u_j)_z$. Conformality gives
$\ip{(u_j)_{zz}}{(u_j)_z}=0$, proving
$\partial_{\bar z}q_j=0$.

The trace formula in Lemma~\ref{lem:c-balance} now reads
\begin{equation}\label{eq:scalar-trace}
 \begin{split}
 q_j(iy)&=-\frac{\lambda_j(y)^2}{2}a_j(s_j(y))
          +i\frac{\lambda_j(y)^2}{2}b(s_j(y)),\\
 g_j(y):=\Ima q_j(iy)&=\frac{\lambda_j(y)^2}{2}b(s_j(y)).
 \end{split}
\end{equation}
The positive factor and the increasing change of parameter preserve
both the number and the order of the signs of $b$.
On the model semicircle, Lemma~\ref{lem:quartic} gives
\[
 z^2(u_j)_{zz}^{\perp}=\frac12(u_j)_{rr}^{\perp}.
\]
Taking its inner product with the real normal $N_j$ shows that
$z^2q_j(z)$ is real there, so $q_j(i)$ is real as well.
Lemma~\ref{lem:corners} also gives
\[
 b(0)=b(\ell)=0.
\]

Integrating the holomorphic form $iq_j(z)(1+z^2)\,dz$ around
$D_+$ gives the moment identity
\begin{equation}\label{eq:scalar-moment}
 \int_{-1}^{1}(1-y^2)g_j(y)\,dy=0.
\end{equation}
To check the boundary terms, on $z=e^{i\theta}$ the form is
\[
 -2\cos\theta\,z^2q_j(z)\,d\theta,
\]
which is real. The positively oriented diameter runs from $i$ to
$-i$. Substituting $z=iy$ and reversing the limits leaves
$\int_{-1}^1(1-y^2)q_j(iy)\,dy$. Its imaginary part must be zero.

Lemma~\ref{lem:endpoint}, with $m=1$, gives
\begin{equation}\label{eq:scalar-endpoint}
 q_j(i)=\frac1\pi\int_{-1}^{1}(1-y^2)R_1(y)g_j(y)\,dy,
 \qquad
 R_1(y)=\frac{1-y+y^2}{(1-y)^2}.
\end{equation}
The function $R_1$ is strictly increasing, since
\[
 R_1'(y)=\frac{1+y}{(1-y)^3}>0
 \qquad(-1<y<1).
\]

\begin{proof}[Proof of Theorem~\ref{thm:torsion}]
If $b$ has one sign, then every $g_j$ has that sign by
\eqref{eq:scalar-trace}. Because the weight in
\eqref{eq:scalar-moment} is positive on $(-1,1)$, the zero moment
forces $g_j\equiv0$. Since $\lambda_j>0$ and $s_j$ maps onto the
whole junction, we obtain $b\equiv0$.

Suppose instead that $b$ changes sign exactly once, from negative
to positive. Choose a separating zero $y_{0,j}$ for each $g_j$.
If an interval of zeros separates the two signs, choose any point
in that interval. Then
\[
 g_j(y)\leq0\quad(y<y_{0,j}),
 \qquad
 g_j(y)\geq0\quad(y>y_{0,j}),
\]
and both signs occur. Subtracting $R_1(y_{0,j})$ from the endpoint
weight is allowed by \eqref{eq:scalar-moment}. Thus
\[
 q_j(i)=\frac1\pi\int_{-1}^{1}(1-y^2)
 \bigl(R_1(y)-R_1(y_{0,j})\bigr)g_j(y)\,dy>0.
\]
The integrand is nonnegative and positive on an open interval, while
Lemma~\ref{lem:endpoint} gives absolute convergence.
At $y=1$, equation~\eqref{eq:scalar-trace} and $b(\ell)=0$ give
\[
 q_j(i)=-\frac{\lambda_j(1)^2}{2}a_j(\ell).
\]
Every $a_j(\ell)$ is therefore negative, contradicting
\eqref{eq:ab-balance}. If $b$ is positive before the transition and
negative afterward, the same integral is negative. Then all
$a_j(\ell)$ are positive, giving the same contradiction.
Since $\nsc(b)\leq1$, these cases prove $b=0$.

Having proved $b\equiv0$, we show that the three faces are planar.
Equation~\eqref{eq:scalar-trace} now makes $q_j$ real on the model
diameter. The function
\[
 H_j(z)=z^2q_j(z)
\]
is real on both boundary arcs, since $(iy)^2=-y^2$ on the diameter.
The maximum principle gives $\Ima H_j=0$ inside the half-disk. The
function $H_j$ is therefore real-valued and holomorphic, hence
constant. Continuity at $0$ gives $H_j\equiv0$ and $q_j\equiv0$.

Write
\[
 h_{xx}=\ip{(u_j)_{xx}}{N_j},\qquad
 h_{xy}=\ip{(u_j)_{xy}}{N_j},\qquad
 h_{yy}=\ip{(u_j)_{yy}}{N_j}.
\]
The equation
\[
 0=q_j=\frac14(h_{xx}-h_{yy}-2ih_{xy})
\]
gives $h_{xx}=h_{yy}$ and $h_{xy}=0$. Harmonicity gives
$h_{xx}+h_{yy}=0$. All three coefficients vanish, so
$\II_j=0$ on each face. Lemma~\ref{lem:planar-triple} completes
the proof.
\end{proof}

\begin{remark}
In $\R^3$, the invariant in Theorem~\ref{thm:oscillation} is
$c_j=a_jb$. In higher codimension, the two second fundamental form
values have additional normal components. Projecting the normal
Hopf field onto a moving unit normal introduces additional
derivatives, so the resulting scalar component is not holomorphic in
general. The all-dimensional criterion therefore uses the full inner
product $c_j$ and the scalar quartic differential.
Theorem~\ref{thm:torsion} is a separate three-dimensional result.
\end{remark}

\section{Admissible variations and the second variation}\label{sec:variations}
\subsection{Normal velocities and their compatibility}
Let $\Sigma$ satisfy Definition~\ref{def:freeboundary}.
An admissible variation is a curve of immersions $X_{j,\varepsilon}$
of the same abstract faces, defined for both signs of $\varepsilon$
near zero. The varied maps agree on the prescribed junction and satisfy
\[
 X_{j,\varepsilon}(S_j)\subset\overline{\B^n},
 \qquad
 X_{j,\varepsilon}(\sigma_j)\subset\Sph^{n-1}.
\]
We take the curves to be twice continuously differentiable in
$\varepsilon$ with values in $C^1$, with $C^{1,\alpha}$ velocities in
the face coordinates. Away from $\varepsilon=0$ the competitors may
lose minimality, conformality, the $120^\circ$ balance, and the free
boundary orthogonality.

Define the velocity by
\[
 U_j=\left.\partial_\varepsilon X_{j,\varepsilon}\right|_{\varepsilon=0}.
\]
Differentiating the common junction trace gives one common velocity
$W$ on the abstract junction. Differentiating the spherical constraint
gives
\[
 U_j|_\Gamma=W,
 \qquad
 \ip{U_j}{X_j}=0\quad\text{on }\sigma_j.
\]
These are variations of the abstract immersed faces. Unrelated
preimages with the same image point can have different velocities.
A single ambient flow is a more restrictive choice, discussed in
Subsection~\ref{sec:ambient-comparison}.

Let $P_j^\perp$ be normal projection on face $j$, and put
$\Phi_j=P_j^\perp U_j$. The normal traces satisfy
\begin{equation}\label{eq:compatibility}
 \Phi_j|_\Gamma=P_j^\perp W\qquad(j=1,2,3)
\end{equation}
for one field $W$. These three normal projections are generally
different. In codimension one, the condition reduces to the scalar
compatibility used in~\cite[Section~2.2]{wang2022curvature}.

\begin{lemma}\label{lem:lift}
Every $C^{1,\alpha}$ normal triple satisfying
\eqref{eq:compatibility} is the normal velocity of an admissible
variation. In particular, normal projections of constant ambient
vectors and of infinitesimal rotations are admissible.
\end{lemma}
\begin{proof}
The map $W\mapsto(P_j^\perp W)_j$ has kernel $\R T$ on the junction,
so compatible normal data have a unique representative $W_0\perp T$.
One convenient formula is
\[
 \mathcal R=T\otimes T+\sum_jP_j^\perp,
 \qquad W_0=\mathcal R^{-1}\sum_j\Phi_j.
\]
If $E$ is the conormal plane, then
$\sum_j\tau_j\otimes\tau_j=\frac32P_E$.
It follows that $\mathcal R$ has eigenvalues $1$ on $\R T$,
$3/2$ on $E$, and $3$ on $(\R T\oplus E)^\perp$.
Its inverse is therefore $C^{1,\alpha}$ along the junction.
For compatible data the displayed formula satisfies
$P_j^\perp W_0=\Phi_j$.

On $\Gamma$ prescribe the tangential correction
\[
 Y_j=W_0-\Phi_j.
\]
At a vertex $p$, Lemma~\ref{lem:corners} gives $p=\pm T$.
Both $W_0$ and $\Phi_j$ are perpendicular to $p$. Thus $Y_j$
lies in the tangent line of the outer edge. Choose a
$C^{1,\alpha}$ tangent field on that edge agreeing with these
two endpoint values, and extend the tangential corrections over
the face. Locally, coefficient data $A(a)$ and $B(b)$ on the
coordinate edges of a quadrant, with common value $v_0$, extend
by $A(a)+B(b)-v_0$. A partition of unity then gives fields
\[
 U_j=\Phi_j+Y_j
\]
that agree on $\Gamma$ and satisfy $\ip{X_j}{U_j}=0$ on $\sigma_j$.

We now construct an admissible variation with this velocity.
Put $\rho_j=1-|X_j|^2$. In an inward transverse coordinate $r$
from an outer edge, including at its vertices, free boundary gives
$\partial_r\rho_j>0$ at $r=0$. In a finite collar we therefore have
$\rho_j\geq cr$ for some $c>0$.
The $C^1$ function $\ip{X_j}{U_j}$ vanishes there at $r=0$, so
its absolute value is at most $C_1r$. Away from this collar,
$\rho_j$ has a positive minimum, and the two estimates combine to
give
\begin{equation}\label{eq:collar-estimate}
 |\ip{X_j}{U_j}|\leq C\rho_j\quad\text{on }S_j.
\end{equation}
Define
\begin{equation}\label{eq:explicit-lift}
 X_{j,\varepsilon}
 =\frac{X_j+\varepsilon U_j}
 {\sqrt{1+\varepsilon^2|U_j|^2}}.
\end{equation}
This has initial velocity $U_j$. Agreement of $X_j$ and $U_j$
on the junction preserves the common trace exactly.
On the outer edge the numerator has squared norm
$1+\varepsilon^2|U_j|^2$, so the image remains on the sphere.
Moreover,
\begin{equation}\label{eq:explicit-deficit}
 \begin{split}
 1-|X_{j,\varepsilon}|^2
 &=\frac{\rho_j-2\varepsilon\ip{X_j}{U_j}}
 {1+\varepsilon^2|U_j|^2}\\
 &\geq\frac{(1-2C|\varepsilon|)\rho_j}
 {1+\varepsilon^2|U_j|^2}.
 \end{split}
\end{equation}
For both signs of sufficiently small $\varepsilon$, this proves
containment in the ball. Uniform $C^1$ convergence to the original
immersions preserves rank two, also at the vertices. Since
\eqref{eq:explicit-lift} is a smooth function of $\varepsilon$ with a
positive denominator for small $\varepsilon$, the curve
$\varepsilon\mapsto X_{j,\varepsilon}$ is smooth with values in
$C^1$; in particular, it has the $C^2$ dependence required in the
definition of an admissible variation.

Finally, $\Phi_j=P_j^\perp v$ is compatible for each constant $v$.
The construction applies even though literal translation can move
the boundary off the sphere. Ambient rotations already preserve
both the sphere and the junction identifications.
\end{proof}

\begin{definition}\label{def:form}
Let $\V$ be the closure in $\bigoplus_jH^1(N\Sigma_j)$ of normal
velocities of admissible variations. The area Hessian defines the
continuous symmetric form $Q$ given in Proposition~\ref{prop:hessian}.
Its Morse index is the largest dimension of a subspace on which $Q$
is negative definite. Its normal kernel and nullity are
\[
 \ker Q=\{\Phi\in\V:Q(\Phi,\Psi)=0\text{ for every }\Psi\in\V\},
 \qquad
 \nul(\Sigma)=\dim\ker Q.
\]
Tangential reparametrizations are not counted, and the condition
$Q(\Phi,\Phi)=0$ by itself does not place $\Phi$ in $\ker Q$.
\end{definition}

Lemma~\ref{lem:lift} describes the compatible regular normal fields,
and $\V$ is defined as their closure rather than through a separate
Sobolev trace characterization at the junction. For the model cone,
Section~\ref{sec:flat} identifies this closure explicitly. The later
index bound uses projected translations, which are already admissible
by the lemma.

\subsection{The area Hessian}
Write $\nabla^\perp_Z\Phi=(D_Z\Phi)^\perp$ for the normal connection.
In a tangent orthonormal frame $(e_1,e_2)$, define
\[
 \Aop_j\Phi=\sum_{a,b=1}^2
 \ip{\II_j(e_a,e_b)}{\Phi}\,\II_j(e_a,e_b).
\]
This endomorphism is independent of the frame. In codimension one,
$\Aop_j\Phi=|\II_j|^2\Phi$. In higher codimension it generally has
more than one eigenspace.
We use the connection Laplacian
\[
 \Delta_j^\perp\Phi
 =\sum_{a=1}^2\left(
 \nabla^\perp_{e_a}\nabla^\perp_{e_a}\Phi
 -\nabla^\perp_{\nabla_{e_a}e_a}\Phi\right),
 \qquad L_j=\Delta_j^\perp+\Aop_j.
\]
Its sign is fixed by Green's formula,
\[
 \int_{\Sigma_j}\ip{\nabla^\perp\Phi}{\nabla^\perp\Psi}\,dA
 =-\int_{\Sigma_j}\ip{\Delta_j^\perp\Phi}{\Psi}\,dA
 +\int_{\partial\Sigma_j}\ip{\nabla^\perp_\nu\Phi}{\Psi}\,ds,
\]
where $\nu$ is the outward conormal on each edge.

\begin{proposition}\label{prop:hessian}
For the variations above, the second variation of area is
\begin{equation}\label{eq:hessian}
 \begin{split}
 Q(\Phi,\Psi)=\sum_{j=1}^3\bigg[&
 \int_{\Sigma_j}\left(
 \ip{\nabla^\perp\Phi_j}{\nabla^\perp\Psi_j}
 -\ip{\Aop_j\Phi_j}{\Psi_j}\right)dA\\
 &-\int_{\sigma_j}\ip{\Phi_j}{\Psi_j}\,ds
 -\int_\Gamma\kappa_j\ip{\Phi_j}{\Psi_j}\,ds\bigg].
 \end{split}
\end{equation}
Here $\kappa_j=\ip{D_TT}{\tau_j}$ uses the outward junction
conormal.
\end{proposition}
\begin{proof}
We calculate first on the diagonal. The interior Jacobi term and
the smooth spherical boundary term are standard, with these
conventions, as in~\cite[Section~3]{fraser2016sharp}.
We give the junction calculation since the three face normals
are not identified in higher codimension.

\proofstep{1}{Differentiate the boundary terms.}
The first variation at the initial surface is
\[
 -\sum_j\int_{\Sigma_j}\ip{H_j}{U_j}\,dA
 +\int_\Gamma\ip{W}{\sum_j\tau_j}\,ds
 +\sum_j\int_{\sigma_j}\ip{U_j}{\eta_j}\,ds=0,
\]
where $H_j$ is the trace mean curvature vector.
The differentiated interior term is
\[
 -\sum_j\int_{\Sigma_j}\ip{L_j\Phi_j}{\Phi_j}\,dA.
\]
On the junction, the acceleration and line-measure terms cancel
by $\sum_j\tau_j=0$. The remaining contribution is
$\sum_j\ip{W}{\dot\tau_j}$.

On an outer edge, the standard boundary calculation gives
\[
 \ip{\Phi_j}{\nabla^\perp_{\eta_j}\Phi_j}-|\Phi_j|^2.
\]
To see this boundary term directly, write $U_j=\Phi_j+\beta_j Z_j$, with $Z_j$
tangent to that edge. The conormal derivative contributes
$\ip{\Phi_j}{\nabla^\perp_{\eta_j}\Phi_j}+\beta_j^2$.
Differentiating $|X_{j,\varepsilon}|^2=1$ twice gives acceleration
flux $-|U_j|^2=-|\Phi_j|^2-\beta_j^2$, so the tangential terms cancel.

\proofstep{2}{Calculate the junction contribution.}
Along $\Gamma$, write
\[
 W=\Phi_j+d_j\tau_j+hT,
 \qquad d_j=\ip{W}{\tau_j},\quad h=\ip{W}{T},
\]
and set $A_j=\II_j(T,T)$, $B_j=\II_j(T,\tau_j)$.
The normal component of the transported conormal derivative is
\[
 (\dot\tau_j)^\perp
 =(D_{\tau_j}U_j)^\perp
 =\nabla^\perp_{\tau_j}\Phi_j-d_jA_j+hB_j.
\]
Here minimality gives $\II_j(\tau_j,\tau_j)=-A_j$.
Differentiating $\ip{\tau_j}{T}=0$ gives its remaining component,
\[
 \ip{\dot\tau_j}{T}
 =-\ip{\tau_j}{D_TW}
 =\ip{\Phi_j}{B_j}-T(d_j)-h\kappa_j.
\]
Its component along $\tau_j$ vanishes because $|\tau_j|=1$.
Thus
\begin{equation}\label{eq:dotconormal}
 \begin{split}
 \ip{W}{\dot\tau_j}={}&
 \ip{\Phi_j}{\nabla^\perp_{\tau_j}\Phi_j}
 -d_j\ip{\Phi_j}{A_j}+2h\ip{\Phi_j}{B_j}\\
 &-hT(d_j)-h^2\kappa_j.
 \end{split}
\end{equation}
We have
\[
 \sum_jd_j=0,\qquad \sum_j\kappa_j=0,
 \qquad B_j=D_T\tau_j+\kappa_jT,\qquad \sum_jB_j=0.
\]
Since each $B_j$ is normal to its face,
$\sum_j\ip{\Phi_j}{B_j}=\ip{W}{\sum_jB_j}=0$.
Also, with $K=D_TT$,
\[
 A_j=K-\kappa_j\tau_j,
 \qquad \ip{\Phi_j}{A_j}=\ip{W}{K}-\kappa_jd_j.
\]
It follows that
\[
 \sum_jd_j\ip{\Phi_j}{A_j}=-\sum_j\kappa_jd_j^2,
 \qquad
 \sum_j\kappa_j|\Phi_j|^2=-\sum_j\kappa_jd_j^2.
\]
The second equality uses
$|\Phi_j|^2=|W|^2-d_j^2-h^2$.
Summing \eqref{eq:dotconormal} now yields
\begin{equation}\label{eq:junction-second}
 \sum_j\ip{W}{\dot\tau_j}
 =\sum_j\left(
 \ip{\Phi_j}{\nabla^\perp_{\tau_j}\Phi_j}
 -\kappa_j|\Phi_j|^2\right).
\end{equation}
\proofstep{3}{Apply Green's formula.}
After integration by parts, the conormal-derivative terms from the
interior Jacobi operator cancel the corresponding boundary terms.
What remains is \eqref{eq:hessian} on the diagonal, and polarization
gives the bilinear formula.

At the stated regularity, these identities can be obtained from
the differentiated area density, which uses first spatial derivatives
of the velocity and acceleration. The Jacobi differentiation and
integration by parts are understood weakly in the interior.
The coefficients and the boundary fluxes above are continuous up
to each edge. Excising the vertices and passing to the limit adds
no flux, since the artificial arcs have vanishing length and the
fluxes are bounded. Equivalently, one can integrate the weak
identity using inner collars. This proves the same formula for
our $C^{2,\alpha}$ immersions and $C^{1,\alpha}$ velocities.
The trace inequality extends it to $\V$.
The codimension-one junction term agrees with
\cite[Theorem~3.2]{wang2022curvature}.
\end{proof}

\begin{remark}\label{rem:sphere-sign}
For the outward normal $n(x)=x$ of the unit sphere, the convention
$h(V,W)=\ip{D_VW}{n}$ gives
\[
 h(V,W)=-\ip{V}{W}.
\]
Tran~\cite[Section~1 and equation~(2.1)]{tran2020index} uses this
convention and writes the support-curvature term with a plus sign.
It is therefore the negative boundary integral in
\eqref{eq:hessian}. The alternative convention
$\widehat h(V,W)=\ip{D_Vn}{W}$ reverses both signs.
Green's formula gives the outward Robin condition
\[
 \nabla^\perp_{\eta_j}\Phi_j=\Phi_j.
\]
For a support sphere of radius $R$, both the boundary coefficient
in the form and the Robin coefficient have magnitude $1/R$.
\end{remark}

The bounded curvature coefficients and the trace inequality
\[
 \|\Phi_j\|_{L^2(\partial S_j)}^2
 \leq\varepsilon\|\nabla^\perp\Phi_j\|_{L^2(S_j)}^2
 +C_\varepsilon\|\Phi_j\|_{L^2(S_j)}^2
\]
show that $Q+C\|\Phi\|_{L^2}^2$ is coercive on $\V$ for large $C$.
The closed form domain contains the smooth normal fields supported
in the face interiors, so it is dense in the corresponding $L^2$
space. Its inclusion into $L^2$ is compact. The associated
self-adjoint operator therefore has compact resolvent, and the
index and nullity are finite.

For a kernel field with sufficient regularity, Green's formula gives
\begin{equation}\label{eq:general-kernel}
 \begin{gathered}
 L_j\Phi_j=0\quad\text{in }\Sigma_j,\\
 \nabla^\perp_{\eta_j}\Phi_j=\Phi_j\quad\text{on }\sigma_j,\\
 \Phi_j|_\Gamma=P_j^\perp W,\qquad
 \sum_j\left(\nabla^\perp_{\tau_j}\Phi_j-\kappa_j\Phi_j\right)=0
 \quad\text{on }\Gamma.
 \end{gathered}
\end{equation}
The interior and spherical equations follow by tests supported away
from the junction. For the last equation, put
$A_j^\Gamma=\nabla^\perp_{\tau_j}\Phi_j-\kappa_j\Phi_j$.
A junction test has normal components $P_j^\perp V$, so its pairing is
\[
 \sum_j\ip{A_j^\Gamma}{P_j^\perp V}
 =\ip{\sum_j A_j^\Gamma}{V}.
\]
On an open subarc, $V$ is arbitrary in $T^\perp$, and the sum on the
right also lies in $T^\perp$; hence the sum vanishes. This derivative
relation belongs to the kernel problem rather than to the form domain.

\section{The spectrum of the flat \texorpdfstring{$Y$}{Y}-cone}\label{sec:flat}
We now specialize to the standard cone. For this calculation, use the
upper half-disk
\[
 \widehat D_+=\{(s,t)\in\R^2:s^2+t^2<1,\ t>0\}.
\]
Here $s$ is the junction coordinate and $t$ is the transverse coordinate.
Write $\widehat\gamma$ for its diameter and $\sigmasp$ for its
semicircle. These coordinates interchange those of the right
half-disk used in the holomorphic arguments.
The face maps are
\[
 x_j(s,t)=se+tv_j,
 \qquad |e|=|v_j|=1,\quad v_j\perp e,\quad \sum_jv_j=0.
\]
The outward junction conormal is $\tau_j=-v_j$. Let
\[
 V=\spanop(YC)\simeq\R^3,
 \qquad F=V^\perp.
\]
Choose an orientation of $V$ and set
\[
 N_j=e\times\tau_j.
\]
Then $\sum_jN_j=0$ and $\ip{N_i}{N_j}=-1/2$ for $i\ne j$.
After choosing an orthonormal basis $E_1,\ldots,E_{n-3}$ of $F$,
every normal field has a unique expression
\begin{equation}\label{eq:flat-normal-frame}
 \Phi_j=f_jN_j+\sum_{\alpha=1}^{n-3}h_{j\alpha}E_\alpha.
\end{equation}
When $n=3$, the sum is absent.
Let $\V_{YC}$ denote the form domain of Definition~\ref{def:form}
for these flat faces.

\begin{lemma}\label{lem:flat-traces}
The space $\V_{YC}$ consists exactly of the triples in
\eqref{eq:flat-normal-frame} whose coefficients belong to
$H^1(\widehat D_+)$ and satisfy
\begin{equation}\label{eq:flat-traces}
 f_1+f_2+f_3=0,
 \qquad h_{1\alpha}=h_{2\alpha}=h_{3\alpha}
 \quad\text{on }\widehat\gamma
\end{equation}
in the trace sense. Notice the difference between the two conditions:
in each additional normal direction the three traces agree, rather
than adding to zero.
\end{lemma}
\begin{proof}
A common junction velocity $W$ gives
\[
 f_j=\ip{W}{N_j},\qquad h_{j\alpha}=\ip{W}{E_\alpha}.
\]
The first traces sum to zero because $\sum_jN_j=0$.
The second traces agree because $E_\alpha$ is the same ambient
normal on every face. Continuity of the Sobolev trace map preserves
these conditions under $H^1$ limits.

Conversely, suppose the trace conditions in \eqref{eq:flat-traces}
hold for a regular triple, and define along the diameter
\[
 W=\frac23\sum_j f_jN_j
       +\sum_\alpha h_{1\alpha}E_\alpha.
\]
For each $k$,
\[
 \ip{W}{N_k}
 =\frac23\left(f_k-\frac12\sum_{j\ne k}f_j\right)=f_k,
 \qquad
 \ip{W}{E_\alpha}=h_{k\alpha}.
\]
These identities give $P_k^\perp W=\Phi_k$ for every $k$.
Lemma~\ref{lem:lift} therefore realizes each regular triple of this
form by an admissible variation.

To pass from regular fields to the full $H^1$ domain, make a constant
orthogonal change in the three sheet indices,
with first coordinate direction $(1,1,1)/\sqrt3$.
In the $f$ sector, the first component is
\[
 \frac{f_1+f_2+f_3}{\sqrt3}.
\]
It has zero diameter trace, and the other two components are
unrestricted. In each $h_\alpha$ sector, the first component is
unrestricted and the other two have zero diameter trace.

An $H^1$ function with zero diameter trace extends oddly across
the diameter to an $H^1$ function on the disk. Extend beyond the
circle, approximate by smooth functions, and take their odd parts.
This gives smooth approximations that still vanish on the diameter.
For an unrestricted function, use even reflection and take even
parts instead. Undoing the orthogonal change gives regular triples
satisfying \eqref{eq:flat-traces}, converging in $H^1$ to the
original triple. The regular-case argument above puts all
approximants in $\V_{YC}$.
\end{proof}

All face and junction curvatures now vanish. Each scalar component
has the form
\[
 q(u,w)=\int_{\widehat D_+}\nabla u\cdot\nabla w\,dA
          -\int_{\sigmasp}uw\,ds.
\]
Define
\[
 \mathcal H_N=H^1(\widehat D_+),
 \qquad
 \mathcal H_D=\{u\in H^1(\widehat D_+):u|_{\widehat\gamma}=0\},
 \qquad
 q_N=q|_{\mathcal H_N},\quad q_D=q|_{\mathcal H_D}.
\]
The letters describe the diameter conditions. The Neumann condition
for $q_N$ appears in the harmonic variational problem, whereas
$q_D$ imposes zero trace there.

The orthogonal change in Lemma~\ref{lem:flat-traces} preserves the
sum of the three gradient integrals and the three outer integrals.
The internal normal sector gives two copies of $q_N$ and one of
$q_D$. Each of the $n-3$ additional directions gives one copy of
$q_N$ and two of $q_D$. Therefore
\begin{equation}\label{eq:form-split}
 Q_{YC}\cong q_N^{\oplus(n-1)}\oplus q_D^{\oplus(2n-5)}.
\end{equation}

\begin{lemma}\label{lem:half-spectrum}
The scalar forms satisfy
\[
 (\Ind q_N,\nul q_N)=(1,1),
 \qquad
 (\Ind q_D,\nul q_D)=(0,1).
\]
For harmonic functions with the corresponding homogeneous Neumann
or Dirichlet condition on the diameter and
$\partial_r u=\varrho u$ on the semicircle, the eigenvalues are
$0,1,2,\ldots$ in the Neumann case and $1,2,3,\ldots$ in the
Dirichlet case. Every eigenvalue in either scalar problem is simple.
\end{lemma}
\begin{proof}
For $u$ in either domain, minimize the Dirichlet energy among
functions in that domain with the same outer trace as $u$. The
admissible affine class is nonempty, and the Poincar\'e inequality
for zero outer trace gives a unique minimizer $h$. Its variational
equation is
\[
 \int_{\widehat D_+}\nabla h\cdot\nabla w\,dA=0
\]
for every allowed $w$ with zero outer trace. Hence $h$ is harmonic,
with zero value on the diameter in the Dirichlet case and the natural
zero conormal derivative there in the Neumann case.
Writing $u=h+w$, we have
\[
 q(h,w)=0,
 \qquad q(w,w)=\int_{\widehat D_+}|\nabla w|^2\,dA>0
 \quad\text{unless }w=0.
\]
The form is therefore an orthogonal sum of its harmonic restriction
and a strictly positive part. This preserves the index. More
explicitly, projection of a negative subspace onto the harmonic
part is injective and remains negative. For a kernel field,
testing against its remainder $w$ gives $q(w,w)=0$, so $w=0$.
The index and kernel are both determined by the harmonic part.

For the harmonic part, use polar coordinates $s=r\cos\theta$,
$t=r\sin\theta$, with $0\leq\theta\leq\pi$.
Even reflection of a Neumann harmonic function across the diameter,
or odd reflection of a Dirichlet one, is weakly harmonic on the full
disk and hence smooth in its interior. The disk Fourier expansion gives
\[
 \begin{array}{ll}
 N:&1,\quad r^k\cos(k\theta),\quad k\geq1,\\
 D:&r^k\sin(k\theta),\quad k\geq1.
 \end{array}
\]
Their outward radial derivatives on $r=1$ give eigenvalues $0,k$
and $k$, respectively. The converse follows from the same Fourier
expansion: an eigenfunction with eigenvalue $\varrho$ can contain
only the coefficient of degree $\varrho$. This also gives simplicity
and completeness, with convergence in the harmonic trace energy norm.

For any one of these modes, the diameter term in Green's formula
vanishes. Hence
\[
 q(u,u)=(\varrho-1)\int_{\sigmasp}u^2\,ds.
\]
Different modes are orthogonal for $q$, as follows from their
boundary orthogonality and Green's formula. The Neumann constant
mode is negative, and its degree-one mode $s$ spans the kernel.
All remaining Neumann modes are positive.
The Dirichlet problem has no negative mode, and its kernel is
spanned by $t$.
\end{proof}

The positive remainder in the proof excludes an additional index
or kernel contribution from functions with zero outer trace.
This point matters in a general Jacobi problem, where a zero
Dirichlet mode can affect the index. See
\cite[Theorems~1.2 and~3.3]{tran2020index}.

\subsection{Kernel fields and rotations}\label{sec:explicit-kernel}
On the diameter, $\partial_\tau=-\partial_t$.
The general junction condition \eqref{eq:general-kernel} becomes
\[
 \sum_j\nabla^\perp_{\tau_j}\Phi_j=0.
\]
The only relation among the $N_j$ is $\sum_jN_j=0$, whereas the
$E_\alpha$ are independent common normal directions. The junction
condition is therefore equivalent to
\begin{equation}\label{eq:flat-derivatives}
 \partial_\tau f_1=\partial_\tau f_2=\partial_\tau f_3,
 \qquad
 \sum_j\partial_\tau h_{j\alpha}=0.
\end{equation}
Lemma~\ref{lem:half-spectrum} and the orthogonal decomposition
show that every kernel component is linear in $(s,t)$.
Write
\[
 f_j=a_js+\beta_jt,
 \qquad h_{j\alpha}=c_{j\alpha}s+d_{j\alpha}t.
\]
The trace conditions give $\sum_ja_j=0$ and
$c_{1\alpha}=c_{2\alpha}=c_{3\alpha}=:c_\alpha$.
Equation~\eqref{eq:flat-derivatives} then gives
$\beta_1=\beta_2=\beta_3=:\beta$ and $\sum_jd_{j\alpha}=0$.
Hence the full kernel is
\begin{equation}\label{eq:kernel-explicit}
 \begin{gathered}
 \Phi_j(s,t)=(a_js+\beta t)N_j
       +\sum_\alpha(c_\alpha s+d_{j\alpha}t)E_\alpha,\\
 \sum_ja_j=0,\qquad \sum_jd_{j\alpha}=0.
 \end{gathered}
\end{equation}
Every field of this form is harmonic and satisfies the trace,
junction derivative, and spherical Robin conditions, so Green's
formula places it in the bilinear kernel. There are three independent
internal parameters and three more for each additional normal
direction.

We next identify these fields as infinitesimal rotations.
For $\omega\in V$, let $\Omega x=\omega\times x$ on $V$ and
extend it by zero on $F$. Since $N_j=-e\times v_j$,
\[
 \ip{\Omega x_j}{N_j}
 =s\ip{\omega}{v_j}-t\ip{\omega}{e}.
\]
Every internal field in \eqref{eq:kernel-explicit} is obtained by
choosing
\[
 \omega=\frac23\sum_j a_jv_j-\beta e.
\]
For the remaining terms, define a linear map $A:V\to F$ by
\[
 Ae=\sum_\alpha c_\alpha E_\alpha,
 \qquad Av_j=\sum_\alpha d_{j\alpha}E_\alpha.
\]
It is well defined because $\sum_jv_j=0$ is the only relation
among the $v_j$, and $\sum_jd_{j\alpha}=0$.
On $V\oplus F$, the map
\[
 \Omega=\begin{pmatrix}0&-A^*\\ A&0\end{pmatrix}
\]
is skew-symmetric and has
\[
 \Omega x_j(s,t)
 =\sum_\alpha(c_\alpha s+d_{j\alpha}t)E_\alpha.
\]
Adding the internal and mixed rotations realizes every field
in \eqref{eq:kernel-explicit}.

\begin{proof}[Proof of Theorem~\ref{thm:flat}]
By \eqref{eq:form-split} and Lemma~\ref{lem:half-spectrum},
\[
 \Ind(YC)=(n-1)\cdot1+(2n-5)\cdot0=n-1,
 \qquad
 \nul(YC)=(n-1)+(2n-5)=3(n-2).
\]
The eigenvalue $0$ occurs only in the $n-1$ Neumann copies.
Every positive integer occurs once in each of the
$(n-1)+(2n-5)$ scalar copies. These are harmonic Steklov
eigenvalues, not the eigenvalues of the $L^2$ Robin Laplacian.
The description \eqref{eq:kernel-explicit} and the rotation
construction above prove the remaining statements.
\end{proof}

For the comparison in Subsection~\ref{sec:ambient-comparison},
a concrete maximal negative subspace will be useful. Its fields are
\begin{equation}\label{eq:negative-constants}
 \Phi_j=\xi_jN_j+Z,
 \qquad \sum_j\xi_j=0,\quad Z\in F.
\end{equation}
The gradients vanish and each semicircle has length $\pi$, so
\[
 Q(\Phi,\Phi)=-\pi\left(\sum_j\xi_j^2+3|Z|^2\right).
\]
This space has dimension $2+(n-3)=n-1$, equal to the index.

\begin{lemma}\label{lem:flat-orthogonal-splitting}
The flat-cone form domain has a $Q$-orthogonal decomposition
\[
 \V_{YC}=\mathcal N\oplus\mathcal K\oplus\mathcal P_+,
\]
where $\mathcal N$ is the negative space in
\eqref{eq:negative-constants}, $\mathcal K=\ker Q$, and $Q$ is
strictly positive on every nonzero element of $\mathcal P_+$.
\end{lemma}
\begin{proof}
Use the scalar decomposition \eqref{eq:form-split}. In each Neumann
copy, Lemma~\ref{lem:half-spectrum} separates the negative constant
mode, the degree-one kernel mode, the higher positive harmonic modes,
and the strictly positive zero-outer-trace remainder. In each
Dirichlet copy, the degree-one mode is the kernel and all higher
harmonic modes and the zero-outer-trace remainder are positive.
Taking the orthogonal sum over all scalar copies gives the stated
decomposition. The negative modes are exactly the constant fields in
\eqref{eq:negative-constants}, and the kernel is the rotation space
described in \eqref{eq:kernel-explicit}.
\end{proof}

\begin{remark}\label{rem:five}
In $\B^3$, the compatible linear triples
$f_j=a_js+\beta_jt$, with $\sum_ja_j=0$, form a five-dimensional
space on which the quadratic value vanishes, although most of these
fields are not in the bilinear kernel. For a compatible constant
triple $\xi$, one has
\[
 Q(f,\xi)=-2\sum_j\beta_j\xi_j,
 \qquad \sum_j\xi_j=0.
\]
Vanishing against every such test requires
$\beta_1=\beta_2=\beta_3$. For example,
\[
 f=(t,0,0),\qquad \xi=(1,-1,0),
 \qquad Q(f,f)=0,\quad Q(f,\xi)=-2.
\]
This distinguishes zero quadratic value from membership in the
bilinear kernel.
\end{remark}

\section{Projected translations and index rigidity}\label{sec:translations}
We return to a general surface satisfying Definition~\ref{def:freeboundary}.
The normal projection $P_j^\perp$ now varies along each face.
The following identity is the counterpart of
\cite[Theorem~3.1]{fraser2016sharp} for the present junction problem.

\begin{theorem}\label{thm:translations}
Let $\Sigma\subset\overline{\B^n}$ be a regular disk-type free boundary
minimal $Y$-surface as in Definition~\ref{def:freeboundary}.
For constant vectors $v,w\in\R^n$, define
$\Phi_j^v=P_j^\perp v$ and $\Phi_j^w=P_j^\perp w$.
Then
\begin{equation}\label{eq:translation-identity}
 Q(\Phi^v,\Phi^v)
 =-2\sum_{j=1}^3\int_{\Sigma_j}|P_j^\perp v|^2\,dA.
\end{equation}
More generally,
\[
 Q(\Phi^v,\Phi^w)
 =-2\sum_{j=1}^3\int_{\Sigma_j}
       \ip{P_j^\perp v}{P_j^\perp w}\,dA.
\]
\end{theorem}
\begin{proof}
Fix $v$ and write $\Phi_j=P_j^\perp v$. Lemma~\ref{lem:lift} makes
these normal data admissible. Standard interior regularity makes each
minimal face smooth away from its boundary, so the normal part of an
ambient translation is a classical Jacobi field in the interior and
\[
 L_j\Phi_j=0
\]
there. This local Jacobi equation does not mean that the translated
boundary remains on the sphere.
Green's formula gives
\[
 \begin{split}
 Q(\Phi,\Phi)={}&\sum_j\int_{\sigma_j}
 \left(\ip{\Phi_j}{\nabla^\perp_{\eta_j}\Phi_j}
                  -|\Phi_j|^2\right)ds\\
 &+\int_\Gamma\sum_j\left(
 \ip{\Phi_j}{\nabla^\perp_{\tau_j}\Phi_j}
                  -\kappa_j|\Phi_j|^2\right)ds.
 \end{split}
\]

\proofstep{1}{Cancellation along the junction.}
Write
\[
 v=\Phi_j+d_j\tau_j+hT,
 \qquad A_j=\II_j(T,T),\quad B_j=\II_j(T,\tau_j).
\]
Differentiating the constant vector $v$ and taking normal components
gives
\[
 \nabla^\perp_{\tau_j}\Phi_j
 =-\II_j(\tau_j,v_j^\top)=d_jA_j-hB_j.
\]
The same geometric identities used in Proposition~\ref{prop:hessian}
give
\[
 \sum_jd_j=0,\quad \sum_j\kappa_j=0,\quad \sum_jB_j=0,
 \qquad \sum_j\ip{\Phi_j}{B_j}=0.
\]
In fact, $\ip{\Phi_j}{B_j}=\ip{v}{B_j}$ because $B_j$ is normal.
With $K=D_TT$, we also have
\[
 \ip{\Phi_j}{A_j}=\ip{v}{K}-\kappa_jd_j,
 \qquad |\Phi_j|^2=|v|^2-d_j^2-h^2.
\]
It follows that
\[
 \sum_jd_j\ip{\Phi_j}{A_j}
 =-\sum_j\kappa_jd_j^2
 =\sum_j\kappa_j|\Phi_j|^2.
\]
Consequently,
\begin{equation}\label{eq:translation-junction}
 \sum_j\left(
 \ip{\Phi_j}{\nabla^\perp_{\tau_j}\Phi_j}
                  -\kappa_j|\Phi_j|^2\right)=0
 \quad\text{on }\Gamma.
\end{equation}
\proofstep{2}{The outer arcs and their endpoints.}
Parametrize one outer arc by arclength with unit tangent $Z$.
Set
\[
 f=\ip{v}{x},\qquad \beta=\ip{v}{Z}.
\]
The vectors $x=\eta_j$ and $Z$ form an orthonormal tangent frame.
Differentiating $\ip{x}{Z}=0$ and $|Z|=1$ gives
\[
 \ip{D_ZZ}{x}=-1,\qquad \ip{D_ZZ}{Z}=0,
 \qquad D_ZZ=-x+\II_j(Z,Z).
\]
Therefore
\[
 f'=\beta,
 \qquad \beta'=-f+\ip{\Phi_j}{\II_j(Z,Z)}.
\]
By minimality and the mixed-curvature identity on the free boundary
\eqref{eq:mixed-fb},
\[
 \nabla^\perp_{\eta_j}\Phi_j
 =-\II_j(\eta_j,v_j^\top)=f\II_j(Z,Z).
\]
Integrating by parts along this arc gives
\[
 \begin{split}
 \int_{\sigma_j}\left(
 \ip{\Phi_j}{\nabla^\perp_{\eta_j}\Phi_j}-|\Phi_j|^2\right)ds
 &=\int_{\sigma_j}(f\beta'+f^2-|\Phi_j|^2)\,ds\\
 &=[f\beta]_{\partial\sigma_j}
   +\int_{\sigma_j}(2f^2-|v|^2)\,ds.
 \end{split}
\]
The last line uses $|v|^2=f^2+\beta^2+|\Phi_j|^2$.

Let $\mu_j(p)$ be the outward unit tangent of the outer boundary
interval at an endpoint $p$. The right-angle sector has inward
outer-arc tangent $-\tau_j(p)$, so $\mu_j(p)=\tau_j(p)$.
At each abstract junction endpoint, the total endpoint contribution
is therefore
\begin{equation}\label{eq:translation-endpoints}
 \ip{v}{p}\ip{v}{\sum_j\mu_j(p)}
 =\ip{v}{p}\ip{v}{\sum_j\tau_j(p)}=0.
\end{equation}
Distinct abstract endpoints are counted separately even if their
images coincide. Combining the two steps gives
\begin{equation}\label{eq:translation-outer}
 Q(\Phi,\Phi)=\sum_j\int_{\sigma_j}
              (2\ip{v}{x}^2-|v|^2)\,ds.
\end{equation}

\proofstep{3}{Convert the remaining flux to an interior integral.}
Consider the ambient vector field
\[
 \mathcal Y_v(x)=|v|^2x-2\ip{v}{x}v.
\]
For an orthonormal tangent frame on a minimal two-dimensional face,
\[
 \diver_{\Sigma_j}\mathcal Y_v
 =\sum_{a=1}^2\ip{|v|^2e_a-2\ip{v}{e_a}v}{e_a}
 =2|v|^2-2|v_j^\top|^2=2|\Phi_j|^2.
\]
The summed junction flux is zero by $\sum_j\tau_j=0$.
On the sphere,
\[
 \ip{\mathcal Y_v}{\eta_j}=|v|^2-2\ip{v}{x}^2.
\]
The divergence theorem and \eqref{eq:translation-outer} prove
\eqref{eq:translation-identity}. Polarization gives the bilinear
version.

The edge fluxes are continuous under our $C^{2,\alpha}$ regularity.
The calculation can first be performed with the vertices excised.
The extra fluxes are bounded times the lengths of the artificial
arcs and tend to zero. The endpoint terms from integration along
$\sigma_j$, in contrast, are retained and canceled in
\eqref{eq:translation-endpoints}.
\end{proof}

\subsection{Proof of the sharp index theorem}\label{sec:rigidity}
\begin{proof}[Proof of Theorem~\ref{thm:rigidity}]
Define
\[
 \mathcal P:\R^n\longrightarrow\V,
 \qquad \mathcal P(v)=(P_1^\perp v,P_2^\perp v,P_3^\perp v).
\]
By Theorem~\ref{thm:translations}, $Q$ is negative definite on
$\operatorname{im}\mathcal P$, because every nonzero normal triple
has positive squared $L^2$ norm. Hence
\[
 \Ind(\Sigma)\geq\operatorname{rank}\mathcal P
              =n-\dim\ker\mathcal P.
\]
A vector in $\ker\mathcal P$ is tangent to every face everywhere:
its continuous normal projections vanish once their $H^1$ classes
do. At any interior junction point $p$,
\[
 \bigcap_{j=1}^3T_p\Sigma_j=\R T(p).
\]
This gives $\dim\ker\mathcal P\leq1$ and therefore
$\Ind(\Sigma)\geq n-1$.

Suppose now that $\Ind(\Sigma)=n-1$. The preceding inequalities give
\[
 \operatorname{rank}\mathcal P=n-1,
 \qquad \dim\ker\mathcal P=1.
\]
Choose a nonzero constant vector $v$ in this kernel. On each face
it is tangent, so for every tangent vector $Z$,
\[
 \II_j(Z,v)=(D_Zv)^\perp=0.
\]
Choose a local orthonormal tangent frame with $e_1=v/|v|$.
Then
\[
 \II_j(e_1,e_1)=\II_j(e_1,e_2)=0.
\]
Minimality gives
\[
 \II_j(e_2,e_2)=-\II_j(e_1,e_1)=0.
\]
All faces are planar. Lemma~\ref{lem:planar-triple} makes them the
full half-unit disks of the standard cone, each covered once.

For the converse, Theorem~\ref{thm:flat} gives $\Ind(YC)=n-1$.
The same argument does not use the value of the index: any nonzero
vector in $\ker\mathcal P$ forces all three faces to be planar.
Hence a nonflat surface has $\ker\mathcal P=\{0\}$, so its projected
translation space has dimension $n$ and its index is at least $n$.
\end{proof}

\subsection{Comparison with ambient variations}\label{sec:ambient-comparison}
There is another natural variation space for an immersed surface.
Let $\V_{\mathrm{amb}}$ be the $H^1$ closure of normal fields induced
by smooth ambient vector fields tangent to the unit sphere.
Their flows preserve the ball, the prescribed junction, and any
additional intersections of the images, so
\[
 \V_{\mathrm{amb}}\subset\V.
\]
Write $\Ind_{\mathrm{amb}}$ for the index of the restricted form and
$\nul_{\mathrm{amb}}$ for the dimension of its bilinear kernel.
The two whole domains are not identified here. We compare only
the quantities needed for the index theorem and the flat model.

\begin{lemma}\label{lem:ambient-translations}
For every constant $v\in\R^n$, the triple
$(P_1^\perp v,P_2^\perp v,P_3^\perp v)$ belongs to
$\V_{\mathrm{amb}}$.
\end{lemma}
\begin{proof}
Set $\rho(x)=1-|x|^2$ and choose a smooth function $\chi$ equal
to one on $(-\infty,1]$ and zero on $[2,\infty)$.
Define
\[
 V_{v,\delta}(x)
 =v-\chi(\rho(x)/\delta)\ip{v}{x}x.
\]
On the sphere its inner product with $x$ is zero.
A smooth ambient cutoff equal to one near the closed ball makes
this field compactly supported without changing its restrictions.
Its flow is an admissible ambient variation.
On face $j$, the difference from the desired normal field is
\[
 E_{j,\delta}
 :=P_j^\perp V_{v,\delta}-P_j^\perp v
 =-\chi(\rho/\delta)\ip{v}{x}x_j^\perp.
\]
It is supported where $\rho\leq2\delta$.
The collar estimate in Lemma~\ref{lem:lift} shows that this set
has area $O(\delta)$, including at the vertices.
Since $x_j^\perp$ vanishes on the outer edge and is $C^1$,
\[
 |x_j^\perp|\leq C\rho,
 \qquad |\nabla^\perp x_j^\perp|\leq C
\]
in that collar. Differentiating the cutoff introduces $1/\delta$
only where $\delta\leq\rho\leq2\delta$. There the factor
$x_j^\perp$ is $O(\delta)$, so
\[
 |E_{j,\delta}|\leq C|v|\delta,
 \qquad |\nabla^\perp E_{j,\delta}|\leq C|v|.
\]
We therefore obtain
\[
 \sum_j\|E_{j,\delta}\|_{H^1(\Sigma_j)}^2
 \leq C|v|^2\delta\longrightarrow0,
\]
which proves the required approximation by ambient normal
velocities, even when the immersed images have additional
intersections.
\end{proof}

\begin{corollary}\label{cor:ambient}
For a surface in Definition~\ref{def:freeboundary},
\[
 \Ind_{\mathrm{amb}}(\Sigma)\geq n-1,
 \qquad
 \Ind_{\mathrm{amb}}(\Sigma)=n-1\quad\Longleftrightarrow\quad\Sigma=YC
\]
up to an orthogonal transformation, with each face covered once.
On the standard cone,
\[
 \Ind_{\mathrm{amb}}(YC)=n-1,
 \qquad \nul_{\mathrm{amb}}(YC)=3(n-2).
\]
\end{corollary}
\begin{proof}
Lemma~\ref{lem:ambient-translations} puts the projected translation
space in $\V_{\mathrm{amb}}$, so Theorem~\ref{thm:translations} and
the same rank argument give the lower bound. Equality again produces
a nonzero constant tangent vector, which forces planarity;
Lemma~\ref{lem:planar-triple} then identifies the surface as $YC$.

For the model, Lemma~\ref{lem:flat-orthogonal-splitting} gives
\[
 \V_{YC}=\mathcal N\oplus\mathcal K\oplus\mathcal P_+,
\]
where $\mathcal N$ is the negative constant space,
$\mathcal K=\ker Q$ is the rotation space, and $Q$ is strictly
positive on every nonzero element of $\mathcal P_+$.
Both $\mathcal N$ and $\mathcal K$ arise from ambient fields.
For the negative space, the field
\[
 U_v(x)=v-\ip{v}{x}x
\]
is tangent to the sphere and has normal component $P_j^\perp v$ on
each flat cone face, since $x$ is tangent to that face. The kernel
consists of actual rotations.
Restriction of a form cannot increase its index, and
$\mathcal N$ attains the full value $n-1$.

If $\Phi$ belongs to the kernel of the restricted form, write
$\Phi=\Phi_-+\Phi_0+\Phi_+$ in the displayed decomposition.
Testing against $\Phi_-\in\mathcal N\subset\V_{\mathrm{amb}}$
first gives $\Phi_-=0$. We then have
\[
 0=Q(\Phi,\Phi)=Q(\Phi_+,\Phi_+),
\]
so positivity on $\mathcal P_+$ gives $\Phi_+=0$ and
$\Phi\in\mathcal K$. Every field in $\mathcal K$ comes from an
ambient rotation, so the restricted kernel is exactly $\mathcal K$.
This gives the stated model nullity as well as the converse index
equality.
\end{proof}

\section{Area and the remaining rigidity question}\label{sec:area}
The geometric assumptions give an area identity independently
of the additional conformal, oscillation, and index hypotheses.

\begin{proposition}\label{prop:area}
For a surface satisfying Definition~\ref{def:freeboundary}, set
\[
 \mathcal A=\sum_j\Area(\Sigma_j),\qquad
 L=\sum_j\Length(X_j(\sigma_j)),\qquad
 \mathcal E=\sum_j\int_{\Sigma_j}|\II_j|^2\,dA,
\]
where areas and lengths are counted on the abstract faces.
Then
\begin{equation}\label{eq:area-curvature}
 2\mathcal A=L=3\pi+\frac12\mathcal E,
 \qquad
 \mathcal A=\frac{3\pi}{2}+\frac14\mathcal E.
\end{equation}
Equality $\mathcal A=3\pi/2$ holds exactly for the standard flat
$Y$-cone, up to an orthogonal transformation.
\end{proposition}
\begin{proof}
We use geodesic curvature with the inward conormal. Along an outer
arc this is $-X_j$, and for a unit tangent $Z$,
\[
 0=Z\ip{X_j}{Z}=1+\ip{X_j}{D_ZZ},
 \qquad \ip{D_ZZ}{-X_j}=1.
\]
Along the junction the inward conormal is $-\tau_j$, so its
geodesic curvature is $-\kappa_j$.
Each disk has two right-angle corners by Lemma~\ref{lem:corners}.
Their exterior angles sum to $\pi$, and Gauss--Bonnet gives
\[
 \int_{\Sigma_j}K_j\,dA+\Length(X_j(\sigma_j))
       -\int_\Gamma\kappa_j\,ds+\pi=2\pi.
\]
The Euclidean Gauss equation for a minimal two-dimensional surface
in any codimension gives
\[
 K_j=-\frac12|\II_j|^2.
\]
Summing and using $\sum_j\kappa_j=0$ yields
$L=3\pi+\mathcal E/2$.

To relate boundary length and area, use
$\diver_{\Sigma_j}X_j^\top=2$, which follows from minimality.
The divergence theorem and free boundary give
\[
 2\Area(\Sigma_j)
 =\Length(X_j(\sigma_j))
       +\int_\Gamma\ip{\gamma_\Gamma}{\tau_j}\,ds.
\]
The junction terms cancel after summation, so $2\mathcal A=L$.

If $\mathcal A=3\pi/2$, then $\mathcal E=0$. Since the squared
norms of the continuous second fundamental forms are nonnegative,
$\II_j=0$ on every face.
Lemma~\ref{lem:planar-triple} gives the standard cone.
Conversely, its three half-unit disks have total area $3\pi/2$.
\end{proof}

\begin{question}\label{question:unrestricted}
Must every disk-type free boundary minimal $Y$-surface in
$\overline{\B^3}$ satisfying Definition~\ref{def:freeboundary}
be the standard flat $Y$-cone, without the matching condition
of Theorem~\ref{thm:main1}, the sign-change conditions of
Theorems~\ref{thm:oscillation} and~\ref{thm:torsion}, or the
index equality of Theorem~\ref{thm:rigidity}?
The same question can be posed in $\B^n$.
\end{question}

Equation~\eqref{eq:area-curvature} gives the lower bound
$\mathcal A\geq3\pi/2$.
Within this geometric class, rigidity would be equivalent to the
reverse universal estimate $\mathcal A\leq3\pi/2$, because
Proposition~\ref{prop:area} characterizes equality. The stationarity
assumptions alone do not provide such an area comparison with the
flat cone.

Any nonflat example in this class would have area strictly greater
than $3\pi/2$. Corollary~\ref{cor:geometric} would force at least
three sign changes in every $c_j$ and at least four in one of them.
In dimension three its common junction torsion would have at least
two sign changes by Theorem~\ref{thm:torsion}.
Its Morse index would be at least $n$ by Theorem~\ref{thm:rigidity}.
These are necessary conditions for a nonflat example, not an
existence result.

\appendix
\section{Conformal coordinates at the junction endpoints}\label{app:corners}
We prove the facewise boundary regularity used in Sections
\ref{sec:oscillation} and~\ref{sec:torsion}.
The geometric assumptions are the $C^{2,\alpha}$ immersion and
edge conditions of Section~\ref{sec:setting}, including rank two
at the vertices. Lemma~\ref{lem:corners} gives the right angles
in the free boundary setting. The argument below applies to one
face at a time and does not identify the three junction parameters.

Two classical regularity results will be used. A positive definite
$C^{1,\alpha}$ metric admits $C^{2,\alpha}$ local harmonic
coordinates, and in dimension two a harmonic conjugate turns them
into isothermal coordinates of the same regularity. See
\cite[Lemma~1.2 and the following discussion]{deturck1981}.
The Kellogg--Warschawski theorem gives $C^{2,\alpha}$ boundary
regularity for conformal maps between Jordan domains with regular
$C^{2,\alpha}$ boundaries. Applying the theorem to the inverse shows
that the boundary derivatives are nonzero.
See~\cite[Chapter~3]{pommerenke1992} and the statement recalled
in~\cite[Introduction]{kalajlamel2020}. We use this smooth-boundary
result only after resolving the corner in an auxiliary planar
domain.

\begin{lemma}\label{lem:corner-straighten}
Suppose that a planar domain near $p$ is bounded by two regular
$C^{2,\alpha}$ arcs meeting at interior angle $\pi/2$.
There is a local conformal coordinate taking $p$ to zero and the
two arcs to the sides of a quadrant. This coordinate and its
inverse are $C^{2,\alpha}$ on the closed sectors, and their
differentials at the vertices are nonsingular.
\end{lemma}
\begin{proof}
\proofstep{1}{Straighten the first arc.}
Extend this arc past $p$ as a regular $C^{2,\alpha}$ curve.
Choose a $C^{2,\alpha}$ Jordan domain whose boundary agrees with
that extension near $p$ and whose interior contains the given
sector locally. The smooth boundary theorem and a M\"obius map
give a conformal map $h$ taking the arc to the real axis, with
\[
 h(p)=0,\qquad h'(p)\ne0.
\]
Both $h$ and its inverse are $C^{2,\alpha}$ through that boundary
point. Choose coordinates so that the sector lies above the real
axis and its first side is the positive real ray. The second side
is a regular $C^{2,\alpha}$ arc perpendicular to the real axis.

\proofstep{2}{Reflect the second arc.}
Parametrize the transformed second arc by arclength as $\beta(s)$,
$s\geq0$, with
\[
 \beta(0)=0,\qquad \beta'(0)=i.
\]
Unit speed gives $\ip{\beta''(0)}{\beta'(0)}=0$ in the real
Euclidean plane. Hence $\beta''(0)$ is real.
Define
\begin{equation}\label{eq:reflected-corner-arc}
 \widehat\beta(s)=
 \begin{cases}
 \beta(s),&s\geq0,\\
 \overline{\beta(-s)},&s\leq0.
 \end{cases}
\end{equation}
The first two derivatives agree at zero, since
\[
 -\overline{\beta'(0)}=\beta'(0)=i,
 \qquad \overline{\beta''(0)}=\beta''(0).
\]
For $s\geq0\geq t$, the difference of the second derivatives is
bounded by
\[
 C(s^\alpha+|t|^\alpha)\leq C'|s-t|^\alpha.
\]
The reflected arc is therefore $C^{2,\alpha}$ through zero. Locally
it is an even graph $x=A(y)$, and the original sector lies on the
side $x>A(y)$.

\proofstep{3}{Straighten the reflected arc symmetrically.}
Close a short portion of this graph into a conjugation-invariant
$C^{2,\alpha}$ Jordan curve. The added joins can be made away from
zero by matching the endpoint derivatives and reflecting the
construction. Let $G$ be the enclosed domain on the right of the
graph. The positive real axis near zero lies inside $G$.
Choose a real point $a_0\in G$ and normalize its Riemann map by
\[
 k:G\longrightarrow\D,
 \qquad k(a_0)=0,\quad k'(a_0)>0.
\]
Uniqueness gives $k(\overline w)=\overline{k(w)}$.
In particular,
\[
 k(0)\in\R\cap\partial\D=\{-1,1\}.
\]
The smooth boundary theorem makes $k$ and $k^{-1}$
$C^{2,\alpha}$ at the boundary with nonzero derivatives.
A real M\"obius map from the disk to the right half-plane taking
$k(0)$ to zero maps the reflected arc to the imaginary axis and
the real interior segment to the real axis. Restricting to the
original sector and composing with $h$ gives the desired quadrant
coordinate, and every differential in this composition is
nonsingular.
\end{proof}

\begin{lemma}\label{lem:equal-corners}
Let a conformal equivalence between planar Jordan domains extend
homeomorphically to their boundaries and take a regular
$C^{2,\alpha}$ right-angle corner to another such corner.
The map and its inverse are $C^{2,\alpha}$ there and have nonzero
complex derivatives. The same conclusion holds in ordinary
$C^{2,\alpha}$ corner coordinates for Riemannian faces with
positive definite $C^{1,\alpha}$ metrics.
\end{lemma}
\begin{proof}
Apply Lemma~\ref{lem:corner-straighten} to the two domains.
In these coordinates the map and its inverse are holomorphic
between quadrant neighborhoods and continuous to their sides.
With corresponding coordinate rays chosen, Schwarz reflection
across the two axes has the formulas
\[
 f(z)=\overline{f(\bar z)},
 \qquad f(z)=-\overline{f(-\bar z)}.
\]
Reflection across the real axis first fills the fourth quadrant,
while reflection across the imaginary axis first fills the second.
Composing the two reflections in either order gives $-f(-z)$ on the
third quadrant. The extensions therefore agree there, and the
reflected map is well defined on a punctured neighborhood of zero.
Boundedness makes the singularity removable.

Applying the same construction to the inverse, the identities
\[
 f^{-1}\circ f=\operatorname{id},
 \qquad f\circ f^{-1}=\operatorname{id}
\]
extend holomorphically to a neighborhood of zero.
Differentiating at zero gives
\[
 (f^{-1})'(0)f'(0)=1.
\]
The derivatives are nonzero. Composing back with the straightening
maps proves the planar statement, including the regularity of
the inverse.

For a Riemannian face, extend its metric coefficients across the
corner and shrink the neighborhood so that positive definiteness
persists. Local isothermal coordinates and their inverses are
$C^{2,\alpha}$. They preserve the intrinsic angle and take the
boundary arcs to regular $C^{2,\alpha}$ planar arcs. The planar
argument then applies.
\end{proof}

\begin{proposition}\label{prop:automatic-regularity}
Let $S$ be a compact $C^{2,\alpha}$ topological disk with boundary
$\partial S=\gamma\cup\sigma$, where the two regular
$C^{2,\alpha}$ arcs meet at two distinct abstract vertices.
Let $X:S\to\R^n$ be $C^{2,\alpha}$ and have rank two everywhere,
including at the vertices, in extendible corner coordinates.
If both intrinsic corner angles are $\pi/2$, there is a conformal
diffeomorphism
\[
 \psi:\overline{D_+}\longrightarrow S,
 \qquad \psi(\gamma_0)=\gamma,\quad\psi(\sigma_0)=\sigma.
\]
The maps $\psi$ and $\psi^{-1}$ are $C^{2,\alpha}$ up to the whole
boundary and have nonsingular differentials. In particular,
$u=X\circ\psi$ is $C^{2,\alpha}$ and has rank two on the closed
half-disk. If $X$ is minimal, then $u$ is harmonic in $D_+$.
\end{proposition}
\begin{proof}
\proofstep{1}{Construct regular conformal coordinates for the metric.}
The induced metric
\[
 g=X^*g_{\R^n}
\]
is positive definite and $C^{1,\alpha}$.
Extend the abstract disk slightly across its edges and corners to
an open topological disk $U$, and extend the metric coefficients in
ordinary half-plane and quadrant charts. After shrinking the charts
if necessary, positive definiteness persists, and a partition of
unity produces a positive definite $C^{1,\alpha}$ metric on $U$.

On a small coordinate neighborhood, choose a $g$-harmonic function
$a\in C^{2,\alpha}$ with $da\ne0$, as provided by local harmonic
coordinates~\cite[Lemma~1.2]{deturck1981}.
Harmonicity gives $d(*_gda)=0$. On a simply connected neighborhood,
choose $b$ with
\[
 db=*_{g}da.
\]
This one-form is $C^{1,\alpha}$, so $b$ is $C^{2,\alpha}$.
Moreover,
\[
 da\wedge db=da\wedge *_gda=|da|_g^2\,dA_g\ne0.
\]
Because the wedge product is nonzero, $(a,b)$ are coordinates. Their
differentials are orthogonal and have the same length, so these
coordinates are isothermal.
The inverse is $C^{2,\alpha}$ as well. With consistent orientation,
the transition maps are holomorphic and define a Riemann-surface
structure on $U$.

\proofstep{2}{Place the compact face in a planar domain.}
By uniformization of the simply connected noncompact Riemann surface
$U$ (see~\cite[Section~16]{mcmullen2019}), there is a biholomorphism
\[
 F:U\longrightarrow\mathcal U,
 \qquad \mathcal U=\C\ \text{or}\ \D.
\]
In the original coordinates, $F$ and $F^{-1}$ are $C^{2,\alpha}$
near the compact face, being holomorphic coordinate changes
composed with the isothermal maps above.
Set
\[
 \Omega=F(\operatorname{int}S),
 \qquad F(S)=\overline\Omega.
\]
The set $\Omega$ is a bounded Jordan domain whose two marked boundary
arcs are $C^{2,\alpha}$ and meet at right angles, since $F$ is
conformal and nonsingular.

\proofstep{3}{Map the half-disk with the prescribed boundary correspondence.}
The Riemann mapping and Jordan boundary extension theorems give
a conformal equivalence
\[
 f:D_+\longrightarrow\Omega
\]
with a homeomorphic extension to the closures. Normalize the
boundary correspondence so that the two corners map to the
prescribed vertices and the diameter maps to $F(\gamma)$.
One may fix the two corners and a third point on the desired arc,
using the three-point normalization of disk automorphisms with
the same cyclic order. Then
\[
 f(\gamma_0)=F(\gamma),
 \qquad f(\sigma_0)=F(\sigma).
\]

At an open boundary point, straighten the corresponding boundary arcs
using auxiliary smooth Jordan domains and the smooth boundary
mapping theorem. Schwarz reflection of the resulting map and
its inverse across the straight intervals gives nonzero boundary
derivatives and $C^{2,\alpha}$ regularity in the original
coordinates. At each corner, Lemma~\ref{lem:equal-corners} gives
that conclusion instead. Compactness combines these local results
on the entire closed boundary.

\proofstep{4}{Return to the original face.}
Define
\[
 \psi=F^{-1}\circ f.
\]
The map has the stated boundary correspondence, and both $\psi$ and
$\psi^{-1}$ are $C^{2,\alpha}$ with nonsingular differentials through
the vertices. By the chain rule, $u=X\circ\psi$ has the same
regularity and remains of rank two. If $X$ is minimal, then the
conformal immersion $u$ is harmonic, which proves the last
conclusion.
\end{proof}

\begin{remark}\label{rem:intrinsic-scope}
Apply Proposition~\ref{prop:automatic-regularity} separately to
the faces of a surface in Definition~\ref{def:freeboundary}, using
the right angles from Lemma~\ref{lem:corners}.
It gives the regular facewise charts of
Definition~\ref{def:boundary-regular}, without identifying their
junction parameters. The reflection is carried out only in an
auxiliary planar coordinate domain.
\end{remark}

\end{document}